\documentclass[11pt,a4paper]{article}
\usepackage[T1]{fontenc}
\usepackage{lmodern}
\usepackage[margin=25mm]{geometry}
\usepackage{amsmath,amssymb,amsthm,mathtools,mathrsfs}
\usepackage{booktabs,array,tabularx}
\usepackage{microtype}
\usepackage{hyperref}
\hypersetup{
  hidelinks,
  bookmarksnumbered=true,
  pdftitle={Sharp spectral comparison and detuning for co-periodic perturbations of deep-water Stokes waves},
  pdfauthor={Bing He and Junfei Cao},
  pdfsubject={Spectral comparison and frequency detuning for deep-water Stokes waves},
  pdfkeywords={Deep-water Stokes waves, co-periodic spectrum, spectral comparison, frequency detuning, Dirichlet--Neumann operator},
  pdflang={en-US}
}
\numberwithin{equation}{section}
\newtheorem{theorem}{Theorem}[section]
\newtheorem{lemma}[theorem]{Lemma}
\newtheorem{proposition}[theorem]{Proposition}
\newtheorem{corollary}[theorem]{Corollary}
\theoremstyle{definition}

\theoremstyle{remark}
\newtheorem{remark}[theorem]{Remark}
\newcommand{\av}[1]{\langle #1\rangle}
\newcommand{\T}{\mathbb T}
\newcommand{\K}{|D|}
\newcommand{\Id}{\mathrm{Id}}
\newcommand{\Lip}{\mathrm{Lip}}
\DeclareMathOperator{\supp}{supp}

\DeclareMathOperator{\meas}{meas}
\allowdisplaybreaks[2]
\title{Sharp spectral comparison and detuning for co-periodic
perturbations of deep-water Stokes waves}
\author{Bing He\thanks{School of Mathematics, Guangdong University of Education,
Guangzhou 510303, Guangdong, China. E-mail: \href{mailto:hzs314@163.com}{hzs314@163.com}}
\and Junfei Cao\thanks{Corresponding author. School of Mathematics,
Guangdong University of Education, Guangzhou 510303, Guangdong, China.
E-mail: \href{mailto:jfcaomath@163.com}{jfcaomath@163.com}. This work was supported by
the National Natural Science Foundation of China (Grant No.~12471176).}}
\date{}
\begin{document}
\maketitle
\begin{abstract}
We study co-periodic perturbations of small-amplitude deep-water
gravity Stokes waves. Replacing the upper-right multiplier $|D_y|$
by $D_y=-i\partial_y$ in conformal coordinates gives a differential
comparison operator. For each nonsquare Fourier index $j\ge2$, the two
physical frequencies $\lambda_{j,\pm}$ continued from $-j\pm\sqrt j$ first
differ from their comparison frequencies at order $e^{2j+2}$, where
$e$ is the first cosine coefficient of the conformal surface profile.
We compute the nonzero leading coefficients explicitly. For nonsquare
$q\ge3$, a flat-water relation detunes as
$\lambda_{4q,+}-\lambda_{q,-}-3\lambda_{q,+}=-6I^3+O_q(I^4)$,
where $I$ is the physical momentum per period length. The case $q=2$
has a different explicit coefficient. At order $I^4$ the gaps divided
by $2\sqrt j$ separate at distinct nonsquare indices $j\ge3$, ruling
out a common three-coefficient dispersion law. These normalized gaps
are strictly ordered over nonsquare $2\le j\le N$ on $0<I\le cN^{-3}$.
\end{abstract}
\smallskip


\section{Introduction}\label{st-introduction}
Frequency combinations that vanish at a flat water surface can become
nonzero along a branch of Stokes waves. Their first nonzero coefficients
give the scale of the resulting frequency denominators. Here we compute
such coefficients for gravity waves in infinite depth and for perturbations
with the same spatial period $2\pi$ as the underlying travelling wave.
At the flat surface the moving-frame frequencies at a positive Fourier
index $j$ are $-j+\tau\sqrt j$, where $\tau=\pm1$ labels the sign of the
square root. We follow the branches at nonsquare indices $j\ge2$, where
these flat eigenvalues are simple and nonzero.

\paragraph{The first contribution from negative Fourier indices.}
After the good-unknown and conformal density changes, the upper-right
entry of the linearization is $|D_y|$, where $y$ is the conformal
coordinate and $D_y=-i\partial_y$. Replacing $|D_y|$ by $D_y$ gives a
differential comparison on the same Fourier space. The two operators
agree on the positive flat blocks. On negative indices their difference
already has order one at zero amplitude. This is the difficulty with
using the comparison to compute the physical frequencies.

Write $e$ for the first cosine coefficient of the conformal surface
profile. Let $\lambda_{j,\tau}$ be the physical frequency and
$\lambda_{\rightarrow,j,\tau}$ its counterpart for the differential
comparison, both continued from the same flat frequency. The arrow
labels the replacement of $|D_y|$ by $D_y$. A mode of signed frequency
$\lambda$ evolves in the moving frame with time factor $e^{-i\lambda t}$.
For each fixed nonsquare $j\ge2$ we prove
\[
 \lambda_{j,\tau}(e)-\lambda_{\rightarrow,j,\tau}(e)
     = C_{j,\tau}e^{2j+2}+O_j(e^{2j+4}).
\]
The error constant and the amplitude interval may depend on $j$.
Theorem~\ref{st-thm:sharp-protection} gives the nonzero coefficient
$C_{j,\tau}$ in closed form, with sign $\tau$.

At Taylor degree $r$, an operator coefficient can shift a Fourier
index by at most $r$. A trip from $j$ to a negative index and back
thus costs at least $2j+2$ degrees. Sharpness requires computing
the coefficient at that degree. The highest Stokes harmonics
give a cancellation in the lower-left entry. The surviving terms
come from the coordinate change itself and from the two-component
inverse at index $-1$. The spectral homotopy in
Section~\ref{st-sec:spectrum} combines the two contributions on
a common domain.

\paragraph{Sixth-order detuning and the limits of a common dispersion law.}
In the constant-transport coordinate, let $v$ be the transport
coefficient and $\kappa$ the lower-left multiplication coefficient.
Two insertions of the leading nonconstant term of $\kappa$ give a
common sixth-order frequency shift for nonsquare $j\ge3$.
With the spatial mean of $\sqrt\kappa$ used as the dispersion
coefficient, the shift through this order equals
$m_s=\operatorname{Var}(\sqrt\kappa)/(2v)$, where the variance uses the
normalized spatial mean in that coordinate. In particular
$m_s=e^6/4+O(e^8)$.

Consider the flat relation
$\lambda_{4q,+}(0)-\lambda_{q,-}(0)-3\lambda_{q,+}(0)=0$,
with $q\ge2$ nonsquare. In the laboratory frame it corresponds to
the five-wave configuration with input wave numbers $q,q,q$ and
output wave numbers $4q,-q$, for which
$3\sqrt q=\sqrt{4q}+\sqrt q$~\cite[Section~2.2]{LucasEtAl21}.
Conjugation puts the minus branch at wave number $-q$ in this
interpretation. The flat quintic coupling coefficient vanishes for
this configuration~\cite[Table~1, case~(iv)]{Lvov97}. We study the
linearized frequency relation along the finite-amplitude Stokes branch.

Write $I$ for the Stokes wave momentum divided by $2\pi$.
The weights $(1,-1,-3)$ cancel the common transport and dispersion
terms, leaving $-3m_s$ through sixth order for $q\ge3$. Since
$I=e^2/2+O(e^4)$, the resulting frequency combination is
\[
 \Delta_q(I)=\lambda_{4q,+}(I)-\lambda_{q,-}(I)-3\lambda_{q,+}(I)
            =-6I^3+O_q(I^4) \qquad(q\ge3).
\]
At $q=2$ the return through the negative sector contributes at the
same order and changes the cubic coefficient to
$-(158+32\sqrt2)/21$.

A law of the form $-j\mathfrak v+\tau\sqrt j\,\mathfrak d+\mathfrak m$
has constant normalized gaps
$s_j=(\lambda_{j,+}-\lambda_{j,-})/(2\sqrt j)$ and affine pair means
$p_j=(\lambda_{j,+}+\lambda_{j,-})/2$ as functions of $j$.
For nonsquare $j\ge3$, the physical branches have this form through
degree six in $e$. At degree eight the normalized gaps differ at
each fixed pair of distinct such indices, and the degree-eight
coefficients of the pair means have strictly decreasing secant
slopes on every increasing triple of such indices.
For instance, the four-branch combination $4\sqrt3(s_3-s_{12})$
has leading term $(17632\sqrt3/20163)I^4$.
Section~\ref{st-sec:consequences} quantifies the best fit when all
three coefficients may be chosen anew at each $I$.
Appendix~\ref{st-sec:polynomial} gives a common interval
$0<I\le cN^{-3}$ for ordering the gaps up to a finite Fourier cutoff
$N$, with $c$ independent of $N$.

\paragraph{Relation to earlier work.}
Longuet-Higgins~\cite[Sections~6--7]{LH78} computed co-periodic
spectra and their leading amplitude-dependent shifts.
Nicholls~\cite[Theorem~5.2]{Nicholls07} proved analytic dependence of selected
nonresonant eigenvalues, and Akers and Nicholls~\cite{AkersNicholls12}
developed expansions for repeated eigenvalues. In the scalar
Burgers--Hilbert equation, Castro, C\'ordoba and
Zheng~\cite[Proposition~3.9 and Corollary~3.10]{CCZ23} used Fourier
support to obtain agreement with their comparison through degree
$2|m|+1$ at Fourier mode $m=n+\operatorname{sign}(n)$ in their
indexing, and computed a nonzero sixth-order error at the lowest
noncentral mode. For water waves the support count gives the
earliest possible order. The highest-harmonic cancellation and
the inverse of the complete two-component block at $-1$ determine
the coefficient at that order for both signs at every nonsquare
index $j\ge2$ of the constrained operator.

For Stokes waves, Berti, Corsi, Maspero and Ventura~\cite{BCMV25}
also use maximal Taylor--Fourier coefficients to determine the
opening of unstable spectral bands at colliding Bloch branches.
Related high-frequency calculations appear in
\cite{CDT22,BMVFirst25}. The rigorous Benjamin--Feir analyses
\cite{NguyenStrauss23,BMV22,BMVFinite23,HsiaoMaspero27} concern
long-wave perturbations and spectrum near zero. Co-periodic
bifurcations at zero and spectral changes toward the highest wave
are studied in \cite{DP25,DMP26}. The results here concern
nonzero simple eigenvalues in the zero Floquet sector, that is,
perturbations with the same period as the Stokes wave.

The water-wave formulation used below is the Hamiltonian surface
formulation of Zakharov~\cite{Zakharov68} and Craig and
Sulem~\cite{CraigSulem93}. We use the analytic Stokes branch
in \cite[Theorem~1.3]{BMVanalytic22} and the full good-unknown
linearization in \cite[(2.9)--(2.12)]{BMV22}. Smoothing of the
Dirichlet--Neumann remainder in one horizontal dimension is already known
\cite{AD15,FG21}. The two-ended Fourier kernel estimate in
Appendix~\ref{st-sec:conformal} supplies the analytic operator
realization needed for the Riesz projections used here~\cite{Kato95}.
The normal form of Feola and Giuliani~\cite{FG21} contains a
sixth-degree Fourier sign term as well as terms of negative spatial order.
Remark~\ref{st-rem:normal-form} explains the identification required
to compare their diagonal entries with these fixed-index frequencies
on the same exact branch. In particular, their decay bound for the
mode-dependent diagonal remainder leaves its sixth-order coefficient
at fixed index undetermined. The direct calculation here gives that
order of detuning in physical momentum. The uniform remainder controls
the competition between the eighth-order gap separation and the
higher-order error as the Fourier cutoff increases.

\paragraph{Organization.}
Section~\ref{st-st-statements} states the spectral laws.
Section~\ref{st-sec:coordinates} constructs the constrained operator
and computes the Stokes coefficients. Section~\ref{st-sec:spectrum}
proves the sharp comparison and the sixth-order detuning.
Section~\ref{st-sec:eighth} treats the eighth-order terms.
Section~\ref{st-sec:consequences} gives the fitting and sublevel
estimates and checks the coefficients in the original conformal
coordinate. The appendices give the analytic kernel estimates and
the index-dependent remainder bounds.

\paragraph{Notation.}
Write $\T=\mathbb R/(2\pi\mathbb Z)$,
$\widehat f(n)=(2\pi)^{-1}\int_0^{2\pi}f(x)e^{-inx}\,dx$,
and $\langle f\rangle=\widehat f(0)$.
We use $P_0f=f-\langle f\rangle$. The multipliers $D$ and $|D|$ have symbols $n$ and $|n|$.
The periodic Hilbert transform $\mathcal H$ has symbol
$-i\operatorname{sign}(n)$, where $\operatorname{sign}(0)=0$. A subscript specifies the spatial
coordinate when needed.
The Sobolev space $H^s=H^s(\T)$ has norm
$\|f\|_{H^s}^2=\sum_n\langle n\rangle^{2s}|\widehat f(n)|^2$,
where $\langle n\rangle=(1+n^2)^{1/2}$ is a Fourier weight.
The zero-mean subspace is $H^s_0$. We write $\|\cdot\|_2$ for
the $H^0=L^2$ norm, summed in squares for vector-valued functions.
Spectral spaces are complexified. Write $\Id$ for the identity,
$T_\vartheta f(x)=f(x+\vartheta)$ for translation, $M_f$ for
multiplication by $f$, and $C_F u=u\circ F$ for composition with a
circle coordinate map $F$.
In a product such as $D\circ b$, $b$ denotes multiplication.
The space of bounded linear maps from $E_1$ to $E_2$ is denoted
$\mathcal B(E_1,E_2)$. A prime on a spatial coefficient denotes
the spatial derivative. Primes on the scalar function $\Delta_q$
refer to $I$. The notation $[e^r]$ takes the coefficient of
degree $r$ in a Taylor series at $e=0$, and $[I^r]$ is used in the
same way for $I$. A subscript in $O_j$ or $O_{j,k}$
allows the bound and neighborhood to depend on the indicated
fixed indices. We write $\meas$ for Lebesgue measure.

\subsection{Physical model and conventions}\label{pc-model}
We consider irrotational gravity waves at infinite depth, with
gravitational acceleration $g=1$ and horizontal period $2\pi$.
Write $U=(\eta,\psi)$, where $\eta(x)$ is the surface height and
$\psi(x)$ is the velocity potential evaluated at the surface.
All potentials below are representatives of
the quotient by spatial constants. At a sufficiently large finite regularity $s$, the real phase space
is $X_s=H^s_0(\mathbb T;\mathbb R)^2$, with coordinates $(\eta,\psi)$.
The zero volume condition is $\langle\eta\rangle=0$. The unnormalized
Dirichlet--Neumann operator $G(\eta)$ is defined by harmonic
extension. Writing $r$ for the vertical coordinate, the potential $\Phi$ in
$\{(x,r):r<\eta(x)\}$ has decaying gradient at infinite depth:
\[
 G(\eta)\psi=(\Phi_r-\eta_x\Phi_x)|_{r=\eta(x)},\qquad
 \Delta\Phi=0,\quad \Phi|_{r=\eta}=\psi.
\]
In particular $G(\eta)1=0$, it is symmetric, and its output has zero
mean. These facts follow from Green's identity, first for smooth finite
energy data and then by continuity. For tangent vectors
$u=(u_\eta,u_\psi)$ and $v=(v_\eta,v_\psi)$,
the Hamiltonian and symplectic form are
\[
 H(U)=\frac12\langle\eta^2+\psi G(\eta)\psi\rangle,
 \qquad \Omega(u,v)=\langle u_\psi v_\eta-u_\eta v_\psi\rangle.
\]
For smooth states, the Hamiltonian vector field $X_H$ is specified
by $\Omega(X_H,h)=-dH[h]$ for every smooth tangent vector $h$.
On a Sobolev scale the evolution is understood at lower regularity
when derivatives are taken.
Its two components are
\begin{equation}\label{pc-equation}
 X_H(U)=\binom{G(\eta)\psi}{P_0f(U)},\qquad
 f(U)=-\eta-\frac12\psi_x^2+
 \frac{(G(\eta)\psi+\eta_x\psi_x)^2}{2(1+\eta_x^2)}.
\end{equation}
The projection chooses the zero-mean representative of the potential.
Integrating the first equation preserves
zero volume. All operations commute with spatial translations and complex
conjugation. The physical momentum per period length is
\begin{equation}\label{pc-momentum}
 \mathcal I(U)=-\langle\psi\eta_x\rangle,
 \qquad d\mathcal I(U)[h]=\Omega(h,-U_x).
\end{equation}
Translation invariance gives $d\mathcal I(U)[X_H(U)]=0$.
Thus total momentum over one period is $2\pi\mathcal I(U)$.
This normalization also applies to the Hamiltonian and symplectic form.

\section{Spectral comparison and detuning}\label{st-st-statements}
\subsection{Exact Stokes branch and signed frequencies}
We use the real analytic small Stokes branch with zero physical mean,
even surface and odd surface potential, supplied by
\cite[Theorem 1.3]{BMVanalytic22}; see also \cite[Theorem 2.1]{BMV22}. For every
finite regularity and positive analytic width needed below, its amplitude
interval may be decreased first. Convergence in the analytic Sobolev
space supplied by that theorem implies convergence in the Fourier--Wiener
norm defined in Appendix~\ref{st-sec:conformal} on any strictly
smaller strip by Cauchy--Schwarz. The conformal
fixed point in Appendix~\ref{st-sec:conformal} is analytic on this
scale. In conformal coordinates write
\[
 X(y)=y+\mathcal H Y(y),\qquad Y(y)=\eta(X(y)),\qquad
 \mathcal H e^{ijy}=-i\,\operatorname{sign}(j)e^{ijy}.
\]
The parameter \(e\) is the \emph{exact} coefficient of \(\cos y\) in \(Y\).
The nonzero first derivative of the Stokes branch makes this a legitimate
local analytic reparameterization. Let \(c_s(e)\) be its wave speed.

Let $U_e=(\eta_e,\psi_e)$ denote the wave at amplitude $e$.
Write $x_{\rm lab}$ for the laboratory horizontal coordinate and $t$
for time. Let \(L_e^{\rm ph}\) be the full linearized evolution
operator in the moving coordinate \(x=x_{\rm lab}-c_s(e)t\).
We write the selected
simple eigenvalues as \(i\beta_{j,\tau}\)
and set
\[
 \lambda_{j,\tau}=-\beta_{j,\tau},\qquad
 \beta_{j,\tau}(0)=j-\tau\sqrt j,\qquad
 j\ge2\ \text{nonsquare},\quad \tau\in\{+1,-1\}.
\]
A corresponding eigenmode evolves as $e^{-i\lambda t}$ in the moving
frame. The index $j$ labels the Fourier mode at $e=0$. The continued
eigenfunction generally has several Fourier components. The parameter
$I=\mathcal I(U_e)$ is the physical momentum per period length, defined in
\eqref{pc-momentum}. In conformal coordinates it is
\(I=-\langle\Psi Y'\rangle_y=\frac12e^2+O(e^4)\), where
\(\Psi=\psi\circ X\) is the conformal surface potential; the identity
and local inversion of $e^2\mapsto I$ are proved in
\eqref{st-eq:action}. The branches are even in $e$. We write
$\lambda_{j,\tau}(I)$ when using this momentum coordinate.

In the transport coordinate $z=F_e(y)$ of \eqref{st-eq:vF}
the transport speed $v(e)$ is constant in space. The lower-left
entry in \eqref{st-eq:Atilde} is multiplication by $\kappa(e,z)$.
With $\langle\cdot\rangle_z$ denoting the normalized spatial mean in
$z$, set
\[
 M=\langle\kappa\rangle_z,\qquad
 d_s=\langle\sqrt\kappa\rangle_z,\qquad
 m_s=\frac{M-d_s^2}{2v}.
\]
The coefficient $v$ is positive for small real $e$ and is distinct
from the Stokes wave speed $c_s$. In $d_s$ we take the positive
square root of the coefficient $\kappa$ near its flat value $1$.

\subsection{Comparison and the common frequency shift}
\begin{theorem}
\label{st-thm:sharp-protection}
Write $\lambda_{\rightarrow,j,\tau}=-\beta_{\rightarrow,j,\tau}$
for the comparison frequency, where $\beta_{\rightarrow,j,\tau}$
is the eigenvalue continued from $\beta_{j,\tau}(0)$. The operator is
$A_{\rightarrow}=\bigl(\begin{smallmatrix}vD_z&D_z\\\kappa&vD_z\end{smallmatrix}\bigr)$
on $L^2(\T;\mathbb C^2)$ with domain $H^1(\T;\mathbb C^2)$.
Here $D_z=-i\partial_z$. In the first
conformal coordinate it replaces the upper-right multiplier $|D_y|$
by $D_y$ and retains every spatial index. The two representations
are related by the same density change as the physical operator,
as shown in Section~\ref{st-sec:spectrum}.
For each nonsquare integer $j\ge2$, put $w=\tau\sqrt j$, $\tau=\pm1$,
and $\gamma_j=-j/(j+1)!$. Then
\begin{equation}\label{st-eq:sharp-protection}
 \lambda_{j,\tau}(e)-\lambda_{\rightarrow,j,\tau}(e)
 = C_{j,\tau}e^{2j+2}+O_j(e^{2j+4}),\qquad
 C_{j,\tau}=\frac{\gamma_j^2}{w}
     \left(1+\frac{2}{(j+1-w)^2-1}\right).
\end{equation}
The coefficient has sign $\tau$. As $j\to\infty$ through nonsquares,
$C_{j,\tau}=\tau j^{3/2}((j+1)!)^{-2}(1+O(j^{-2}))$.
In momentum variables the difference is
$2^{j+1}C_{j,\tau}I^{j+1}+O_j(I^{j+2})$.
\end{theorem}

The factor $\gamma_j$ is the first nonzero Taylor coefficient for
the crossing from $j$ to $-1$ under the transport coordinate map
(Lemma~\ref{st-lem:highest}).
Appendix~\ref{st-sec:polynomial} gives absolute bounds on discs of
radius proportional to $\operatorname{dist}(2\sqrt j,\mathbb Z)^{1/2}j^{-1/4}$,
where $\operatorname{dist}$ is ordinary distance.
The restriction to nonsquares avoids a double eigenvalue at the
flat wave. If $j=m^2$ with $m\ge2$, its branch with label $\tau$
coincides at zero amplitude with the branch of the opposite sign
at $(m-\tau)^2$.

\begin{proposition}\label{st-prop:individual}
For each fixed nonsquare integer $j\ge3$ and $\tau\in\{+1,-1\}$,
\[
 \lambda_{j,\tau}(e)
  =-jv(e)+\tau\sqrt{jM(e)}+\frac{1-\tau\sqrt j}{4}e^6+O_j(e^8).
\]
Equivalently, to this order,
\begin{equation}\label{st-eq:variance-expansion}
 \lambda_{j,\tau}(e)=-jv(e)+\tau\sqrt j\,d_s(e)+m_s(e)+O_j(e^8),
 \qquad m_s(e)=\tfrac14 e^6+O(e^8).
\end{equation}
\end{proposition}

The exact identity $2vm_s=\langle(\sqrt\kappa-d_s)^2\rangle_z$ gives
$m_s>0$ for all sufficiently small nonzero real $e$. Indeed the first
Fourier coefficient of $\kappa$ is $e^3+O(e^4)$ by \eqref{st-eq:kap3}.
Section~\ref{st-sec:return} identifies this variance with the common
offset through sixth order.

At $j=2$, the comparison error has exactly the order of the two
cubic insertions in the proof of Proposition~\ref{st-prop:individual}.
Explicitly,
\begin{equation}\label{st-eq:boundary}
 \lambda_{2,\tau}-\lambda_{\rightarrow,2,\tau}
 =\frac{1+2\tau\sqrt2}{21}e^6+O(e^8).
\end{equation}
For $j\ge3$ the comparison first changes the frequency at order eight
or higher in $e$.

\subsection{Detuning and normalized gaps}
\begin{theorem}[Sixth-order detuning]\label{st-thm:detuning}
Fix a nonsquare integer \(q\geq2\). There is \(I_q>0\) such that the three
indicated simple branches of the full, volume-constrained water-wave
linearization on the potential quotient satisfy, for \(0\leq I<I_q\),
\begin{equation}\label{st-eq:main}
 \Delta_q(I):=\lambda_{4q,+}(I)-\lambda_{q,-}(I)-3\lambda_{q,+}(I)
          =-c_q I^3+O_q(I^4),
\end{equation}
where $c_2=(158+32\sqrt2)/21$ and $c_q=6$ for $q\ge3$.
The remainder is analytic in \(I\). For $q\ge3$,
\[
 \Delta_q'(I)=-18I^2+O_q(I^3),
 \qquad -9I^3\leq\Delta_q(I)\leq-3I^3
\]
after a reduction of \(I_q\). The radius and remainder constant may depend on $q$.
\end{theorem}

\begin{theorem}[Eighth-order separation of normalized gaps]
\label{st-thm:gap-variation}
For each nonsquare $j\ge2$, write
$s_j(I)=(\lambda_{j,+}(I)-\lambda_{j,-}(I))/(2\sqrt j)$
for the normalized branch gap.
For any fixed distinct nonsquare $j,k\ge3$,
\begin{equation}\label{st-eq:gap-variation}
 s_j(I)-s_k(I)
 =16\bigl(\mathfrak f(j)-\mathfrak f(k)\bigr)I^4+O_{j,k}(I^5),
 \qquad
 \mathfrak f(j)=f(j)+\frac7{1056}\mathbf 1_{\{3\}}(j),
\end{equation}
where
$f(x)=1/(8x)-9/[64(4x-1)]-3/[64(4x-9)]$ and
$\mathbf 1_{\{3\}}$ is one at $j=3$ and zero elsewhere.
The added term at $j=3$ in $\mathfrak f$ is the correction from
negative Fourier indices. The function $\mathfrak f$ is strictly
decreasing on nonsquare indices $j\ge3$. In particular $s_j(I)>s_k(I)$ for $j<k$ and sufficiently
small positive $I$, with a neighborhood depending on the pair.
\end{theorem}

\begin{corollary}
\label{st-cor:eighth-detuning}
For the four branches with indices $3$ and $12$ and both signs, set
\[
 Q(I)=2\lambda_{3,+}-2\lambda_{3,-}-\lambda_{12,+}+\lambda_{12,-}
 \qquad \text{with all frequencies evaluated at }I .
\]
Then as $I\downarrow0$,
\begin{equation}\label{st-eq:eighth}
 Q(I)=\left(\frac{3026}{6721}+\frac{14}{33}\right)\sqrt3\,I^4
        +O(I^5)=\frac{17632\sqrt3}{20163}I^4+O(I^5).
\end{equation}
The first summand is the contribution of the differential comparison
operator and the second is the change caused by the negative
spatial indices.
\end{corollary}

\section{Constrained linearization and Stokes coefficients}\label{st-sec:coordinates}
\subsection{The constrained operator}\label{st-sec:physical}
\subsubsection*{Good unknown and conformal density}
Throughout this section the background is a fixed exact Stokes wave.
The physical moving coordinate is $x=x_{\rm lab}-c_s(e)t$.
We use $y$ for the conformal boundary coordinate, with $x=X_e(y)$,
and later $z=F_e(y)$ to make the transport coefficient constant.
Both spatial maps are static at the fixed amplitude $e$.
For a physical perturbation $(h,k)$, let $B_{\rm ph}$ and
$V_{\rm ph}$ be the vertical and horizontal fluid velocities on
the surface, respectively. They are given by
\[
 B_{\rm ph}=\frac{G(\eta)\psi+\eta_x\psi_x}{1+\eta_x^2},\qquad
 V_{\rm ph}=\psi_x-B_{\rm ph}\eta_x .
\]
The full shape derivative, or the full linearization in
\cite[(2.9)--(2.11)]{BMV22}, gives, after the good unknown
$(h,p_{\rm g})$, where $p_{\rm g}=k-B_{\rm ph}h$, the operator
\[
 L_{\rm good}=
 \begin{pmatrix}
   \partial_x\circ(c_s-V_{\rm ph})&G(\eta)\\
   -[1+(V_{\rm ph}-c_s)B_{{\rm ph},x}]&
                      (c_s-V_{\rm ph})\partial_x
 \end{pmatrix}.
\]
The second component is read modulo constants. The \(P_0\) in
\eqref{pc-equation} selects a representative of this quotient.

Use the conformal boundary map \(X(y)=y+\mathcal H Y(y)\) and
define \(\mathcal T_e(h,k)=(\xi,p)\) by
\[
 \xi=X'(h\circ X),\qquad
 p=p_{\rm g}\circ X .
\]
This map is linear in the perturbation at the fixed background.
Its inverse is
\[
 h=(\xi/X')\circ X^{-1},\qquad
 k=(p+B\xi/X')\circ X^{-1},
 \quad B=B_{\rm ph}\circ X,
\]
with a final choice of representative if desired.
Change of variables gives the mean and symplectic identities
below, with $\delta$ denoting the differential in perturbation
coordinates.
\[
 \av{\xi}_y=\av{h}_x,\qquad
 \av{\delta p\wedge\delta \xi}_y
      =\av{\delta k\wedge\delta h}_x.
\]
Adding a physical potential constant adds precisely the same constant
to \(p\). Thus the constrained physical space is carried exactly to
\[
 \{\xi\in H^s:\av{\xi}=0\}\ \times\ (H^s/\mathbb C).
\]
The zero-mean representative of \(p\) need not correspond to the
zero-mean representative of \(k\).

Let $Z(\zeta)$ be the conformal map from the lower half-plane to the
fluid domain, with boundary value $Z(y)=X(y)+iY(y)$; $\zeta$ denotes
the complex conformal coordinate. The steady kinematic condition
makes the relative streamfunction constant at the surface. After
subtracting the uniform far-field potential $-c_s\zeta$, the periodic
imaginary part is bounded and has constant boundary data. The harmonic
extension in the lower half-plane is that constant. Up to additive
constants the relative holomorphic potential is therefore $-c_s\zeta$.
Write $J=X'^2+Y'^2$, which equals $|Z'|^2$ for real amplitudes.
The polynomial expression defines $J$ at complex amplitudes.
For the real wave it follows that
\begin{equation}\label{st-eq:velocities}
 (V_{\rm ph}-c_s)\circ X=-\frac{c_sX'}J,\qquad
 B=-\frac{c_sY'}J .
\end{equation}
The energy identity \eqref{st-eq:energy} gives
$X'\,(G(\eta)f)\circ X=\K(f\circ X)$. The transport entries
follow from $(c_s-V_{\rm ph})\circ X=(c_s/J)X'$ and
$\partial_y(f\circ X)=X'(\partial_x f)\circ X$.
On the constrained quotient, $\mathcal T_e$ gives $L_e^{\rm c}=\mathcal T_e L_e^{\rm ph}\mathcal T_e^{-1}$.
Its formula, extended to the ambient two-component space, is
\[
 L_e^{\rm c}=\begin{pmatrix}\partial_y\circ b&\K\\-a&b\partial_y\end{pmatrix},
 \qquad
 b=\frac{c_s}J,\qquad a=\frac{1-c_sB_y/J}{X'}.
\]
Let \(C=\operatorname{diag}(1,i)\). The convenient spectral operator is
\begin{equation}\label{st-eq:A}
 A_e=-iC^{-1}L_e^{\rm c}C
       =\begin{pmatrix}D\circ b&\K\\a&bD\end{pmatrix}.
\end{equation}
The matrix \(C\) is a complex similarity used to pass from the
evolution eigenvalue \(i\beta\) of \(L_e^{\rm c}\) to the eigenvalue
\(\beta\) of \(A_e\).

\subsubsection*{A common closed domain}
The ambient complex Hilbert space for \eqref{st-eq:A} is
\(\mathcal X=L^2(\T;\mathbb C^2)\), with common domain
\(\mathcal D=H^1(\T;\mathbb C^2)\). At zero amplitude the Fourier
block at index $n$ is
\[
 A_0(n)=\begin{pmatrix}n&|n|\\1&n\end{pmatrix}.
\]
For any fixed nonreal \(\zeta\), its inverse at large \(|n|\) is
\[
 (A_0(n)-\zeta)^{-1}
 =\frac1{(n-\zeta)^2-|n|}
       \begin{pmatrix}n-\zeta&-|n|\\-1&n-\zeta\end{pmatrix}
       =O(|n|^{-1})
\]
as a two-component matrix. This proves the graph estimate
\[
 \|u\|_{H^1}\leq C(\|A_0u\|_2+\|u\|_2).
\]
The coefficients $a,b$ give an analytic family
\(A_e-A_0=O(e)\) in \(\mathcal B(H^1,L^2)\).
Absorbing its relative bound gives the same estimate for \(A_e\).
A graph-convergent sequence is then Cauchy in \(H^1\), proving
closedness on the stated common domain. A Neumann argument at one
fixed nonreal \(\zeta\) gives a compact resolvent.
The same argument on any finite list of higher Sobolev levels holds
after those levels are fixed and the amplitude radius is reduced.

This is also the distributional maximal domain. If
$(\xi,p)\in L^2$ and $A_e(\xi,p)=f\in L^2$ in distributions, the
second row gives $Dp=b^{-1}(f_2-a\xi)\in L^2$. The first row then gives
$D(b\xi)=f_1-|D|p\in L^2$. Since $b$ is smooth and nonvanishing,
both $p$ and $\xi$ belong to $H^1$. For the reduced operator below
the omitted constant in the second row is the bounded functional
$\langle a\xi-(Db)p\rangle$, so the same argument applies to the
zero-mean representatives. The real coordinate maps preserve $H^1$
and their identities extend to distributions by duality.

\subsubsection*{The full mean block and the potential quotient}
Set \(\mathcal X_0=(L^2_0)^2\), \(e_\xi=(1,0)\), \(e_p=(0,1)\).
Here $e_\xi$ spans the surface-mean direction excluded by the volume
constraint, and $e_p$ spans the constant-potential direction identified
with zero in the quotient. In the fixed decomposition
\(\mathbb C e_\xi\oplus\mathcal X_0\oplus\mathbb C e_p\),
the operator has the exact decomposition
\begin{equation}\label{st-eq:mean-block}
 A_e=\begin{pmatrix}
       0&0&0\\ b_0&A_{\rm red}&0\\ \bar a&\varrho&0
     \end{pmatrix},
 \quad b_0=(Db,P_0a),\quad\bar a=\av a,
\end{equation}
where
\begin{align*}
 A_{\rm red}(\xi,p)&=(D(b\xi)+\K p,\ P_0(a\xi+bDp)), \\
 \varrho(\xi,p)&=\av{a\xi-(Db)p}.
\end{align*}
Integration by parts gives the formula for \(\varrho\). This row is
\(L^2\)-bounded, while \(b_0\) is smooth. The domain is
\(\mathbb C\oplus(H^1\cap\mathcal X_0)\oplus\mathbb C\).
The reduced operator is closed: convergence of its graph also gives
convergence of the full graph after adjoining its bounded last row.

\begin{proposition}\label{st-prop:quotient}
The full and reduced operators in \eqref{st-eq:mean-block} have the same
nonzero spectrum. At isolated finite spectral sets away from zero
their Jordan structures and algebraic multiplicities agree.
\end{proposition}
\begin{proof}
For \(\zeta\ne0\), if \(R_{\rm red}(\zeta)=(A_{\rm red}-\zeta)^{-1}\) exists, the
unique solution of \((A_e-\zeta)(h_0,u,p_0)=(f_\xi,f_0,f_p)\) is
\[
 h_0=-f_\xi/\zeta,\quad
 u=R_{\rm red}(\zeta)(f_0+b_0f_\xi/\zeta),\quad
 p_0=(\varrho(u)+\bar a h_0-f_p)/\zeta .
\]
Every product is legitimate on the specified strong and weak spaces.
Conversely, full surjectivity with forcing \((0,f_0,0)\) gives reduced
surjectivity. A reduced kernel vector \(u\) lifts to
\((0,u,\varrho(u)/\zeta)\), so full injectivity gives reduced injectivity.
The closed graph theorem supplies the reduced bounded inverse.

On a finite reduced generalized spectral space \(E\) away from zero,
\(A_{\rm red}|_E\) is invertible and
\[
 Su=(0,u,\varrho A_{\rm red}^{-1}u),\qquad A_eS=SA_{\rm red}.
\]
For a full generalized spectral vector at \(\zeta\ne0\), its first
constant component is zero: the first row of \((A_e-\zeta)^m\) is
\((-\zeta)^m\) on that component. Projection onto the middle component
is therefore injective on the corresponding spectral space; its
kernel could only be a potential constant, which belongs to the
zero-spectral sector. The displayed intertwiner gives
the reverse inclusion and surjectivity. Hence the finite Jordan
structures coincide. The intertwiner retains the finite-rank
coupling \(\varrho\) to the constant potential direction.
\end{proof}

Conserved quantities also lead to volume and momentum constraints in
other conformal co-periodic formulations; see~\cite[Remark~6]{DP25}.
The triangular resolvent above gives the nonzero-spectrum identification
for the variables and sign convention used here.

\subsubsection*{Simple branches, reality, and parity}
For nonsquare \(j\geq2\), \(\beta_0=j-\tau\sqrt j>0\).
Negative spatial blocks have eigenvalues \(n\pm\sqrt{|n|}\leq0\).
If a positive block \(n\) has the same eigenvalue, write
\(j=dr^2\) with \(d>1\) squarefree and compare rational and
squarefree-root parts in \(n\pm\sqrt n=j-\tau\sqrt j\).
It follows that \(n=j\) and the same half is selected.
Thus \(\beta_0\) is simple in the full operator.
For comparison, if $j=m^2$ with $m\ge2$, then
$j-\sqrt j=(m-1)^2+(m-1)$ and
$j+\sqrt j=(m+1)^2-(m+1)$.
The nonsquare condition excludes these genuine flat degeneracies.

A fixed small resolvent circle and the common-domain Neumann argument
give a rank-one analytic Riesz branch. Higher Sobolev levels needed
for a finite number of derivatives are fixed before decreasing its
radius. This constructs smooth eigenvectors with finite parameter
jets without asserting operator-norm analyticity of composition on
one Sobolev level.

For real \(e\), the coefficients \(a,b\) are real-valued and even in \(y\). The
antilinear involution \(\mathcal Rf(y)=\overline{f(-y)}\) commutes with
\(D,\K\), and \(A_e\). Uniqueness of the simple branch implies
\(\beta(e)\in\mathbb R\). Separately, uniqueness of the Stokes branch
and a spatial \(\pi\)-translation give
\[
 Y_{-e}(y)=Y_e(y+\pi),\qquad
 A_{-e}=T_\pi^{-1}A_eT_\pi .
\]
Therefore \(\beta(-e)=\beta(e)\).

\subsection{Conformal coefficients and constant transport}\label{st-sec:stokes}
\subsubsection*{Bernoulli's constant is fixed by physical volume}
Let \(\Psi=\psi\circ X\). The relative complex potential $\mathcal W_{\rm rel}$ satisfies
$\mathcal W_{\rm rel}(Z(\zeta))=-c_s\zeta$; its boundary real part gives
\[
 \Psi=c_s(X-y)=c_s\mathcal H Y
       \quad\hbox{modulo constants}.
\]
First retain an unspecified Bernoulli constant \(\mathcal B\):
\[
 (\mathcal B-2Y)|Z'|^2=c_s^2,\qquad
 Z'=1+\K Y+iY'.
\]
Here \(Z'\) has nonpositive spatial frequencies and no zeros for
small waves. Thus
\((\mathcal B-2Y)Z'=c_s^2/\overline{Z'}\) has nonnegative frequencies
and constant coefficient \(c_s^2\).
Its imaginary part is \((\mathcal B Y-Y^2)'\).
For a function with nonnegative Fourier support, the real part is
its real constant term minus \(\mathcal H\) applied to its imaginary part. Since
\(\mathcal H\partial_y=\K\), comparison of real parts gives
\[
 (\mathcal B-c_s^2)+2\mathcal B\K Y-2Y-2Y\K Y-\K(Y^2)=0.
\]
The physical zero-volume condition is
\[
 0=\av{\eta}_x=\av{YX'}_y=\av{Y+Y\K Y}.
\]
Taking the mean in the real-part identity therefore gives
\(\mathcal B=c_s^2\). In particular the exact Babenko identity is
\begin{equation}\label{st-eq:Babenko}
 (c_s^2\K-1)Y=Y\K Y+\tfrac12\K(Y^2).
\end{equation}

\subsubsection*{Taylor--Fourier filtration}
Write
\[
 Y=\sum_{n\geq1}e^nY_n,\quad
 c_s^2=1+\sum_{r\geq1}\sigma_re^r,\quad
 Y_1=\cos y,\quad \widehat Y_n(\pm1)=0\ (n>1).
\]
The last normalization refers to the coefficients of the first
cosine harmonic. Comparing coefficients in \eqref{st-eq:Babenko} gives
\begin{equation}\label{st-eq:recursion}
 (\K-1)Y_n
 =\sum_{a+b=n}\left(Y_a\K Y_b+\tfrac12\K(Y_aY_b)\right)
        -\sum_{r=1}^{n-1}\sigma_r\K Y_{n-r}.
\end{equation}
The half-period symmetry of the Stokes branch makes $c_s$ even, so
$\sigma_1=0$. The first cosine coefficient determines the new \(\sigma_{n-1}\); the
remaining Fourier coefficients determine \(Y_n\). Induction gives
\[
 \supp\widehat Y_n\subset\{|k|\leq n\}.
\]
Indeed products add bandwidths, \(\K\) and derivatives do not enlarge
them, and the sole even kernel of \(\K-1\) is fixed by the normalization.
This is also a special case of the Stokes support structure in
\cite{BCMV25}.

\begin{lemma}\label{st-lem:filtration}
Every Taylor coefficient of order \(r\) of the coefficient
multiplications, the near-identity circle maps used below, their
inverse maps, and their densities has spatial bandwidth at most \(r\).
The order-\(r\) coefficient of a composition operator moves a Fourier
label by at most \(r\).
\end{lemma}
\begin{proof}
The assertion follows first for \(X',Y',J,a,b\) from the bandwidth bound
for \(Y_n\).
It is preserved by products and by inverses of scalar series whose
constant coefficient is a nonzero scalar. Spatial integration with
zero mean does not enlarge support.
For a near-identity map \(F=y+f\), with \(e_j(y)=e^{ijy}\), expansion of
\[
 C_Fe_j=e^{ijy}\sum_{k\geq0}\frac{(ij f)^k}{k!}
\]
shows that each coefficient of total order \(r\) has offset at most
\(r\). The inverse map follows inductively from
\(F\circ F^{-1}=\Id\); each term has the same total-order filtration.
The Jacobian and inverse Jacobian have the scalar-series property.
Only finite jets of these formulas are used on Sobolev scales,
with enough spatial derivatives fixed in advance.
\end{proof}

\subsubsection*{Third-order table}
Put \(C_n=\cos(ny)\), \(S_n=\sin(ny)\). At second order the right-hand
side of \eqref{st-eq:recursion} is \(1/2+C_2\).
At third order it is \(C_1+3C_3-\sigma_2C_1\).
The first-harmonic solvability condition therefore yields
\[
 Y_2=-\tfrac12+C_2,\qquad \sigma_2=1,\qquad Y_3=\tfrac32 C_3.
\]
Using \eqref{st-eq:velocities} and \eqref{st-eq:A} gives the full table
\begin{align}
 Y&=eC_1+e^2(-\tfrac12+C_2)+\tfrac32e^3C_3+O(e^4),\label{st-eq:Y3}\\
 c_s&=1+\tfrac12e^2+O(e^4),\nonumber\\
 J&=1+2eC_1+e^2(1+4C_2)+e^3(4C_1+9C_3)+O(e^4),\nonumber\\
 b&=1-2eC_1+e^2(\tfrac32-2C_2)+e^3(C_1-3C_3)+O(e^4),\nonumber\\
 B&=eS_1+e^2S_2+e^3(\tfrac12S_1+\tfrac32S_3)+O(e^4),\nonumber\\
 bB_y&=eC_1+e^2(C_2-1)+e^3(-C_1+\tfrac32C_3)+O(e^4),\nonumber\\
 a&=1-2eC_1+e^2(2-2C_2)+e^3(2C_1-3C_3)+O(e^4).\nonumber
\end{align}
Each row is obtained by finite division and trigonometric products,
with \(a=(1-bB_y)/(1+\K Y)\). The remainders refer to the actual
analytic Stokes branch in any preselected finite strong space.

\subsubsection*{Constant transport by a second density change}
Set
\begin{equation}\label{st-eq:vF}
 v(e)=\av{1/b}^{-1}>0,\qquad F'(y)=v/b(y),\qquad
 \av{F(y)-y}=0 .
\end{equation}
The average condition on \(F'\) gives a degree-one circle map.
This \(F:y\mapsto z\) is distinct from the physical conformal
map \(X:y\mapsto x\).
On the real slice the spatial $L^2$ adjoint of $C_F$ is
$C_F^\sharp=C_F^{-1}M_{1/F'}$.
At complex parameters $\sharp$ denotes this ordered continuation,
not a Hilbert adjoint; its analytic justification is given in
Appendix~\ref{st-sec:conformal}.
For old variables \((\xi,p)\) and new variables
\((\widetilde \xi,\widetilde p)\), define
\begin{equation}\label{st-eq:density}
 \xi=F'C_F\widetilde \xi,\qquad p=C_F\widetilde p .
\end{equation}
This density map is invertible and symplectic:
change of variables preserves both \(\av{\xi}\) and the symplectic form.
The potential quotient is implemented by \(P_0C_F\) and
\(P_0C_F^{-1}\) as in Section~\ref{st-sec:conformal}.
For fixed small real \(e\), the maps are bounded in both directions
on \(L^2\) and \(H^1\), and carry the corresponding closed domains.

With $bF'=v$, the two diagonal entries in the new coordinate are
$vD_z$. The lower-left entry is multiplication by $(aF')\circ F^{-1}$,
and the upper-right entry is the density congruence of $|D|$.
The full transformed operator is
\begin{equation}\label{st-eq:Atilde}
 \widetilde A_e
 =\begin{pmatrix}
       vD_z&D_z+\mathcal H_e\\
       \kappa(z)&vD_z
   \end{pmatrix},
 \quad
 \kappa(F(y))=\frac{v a(y)}{b(y)},\quad
 \mathcal H_e=C_F^\sharp(|D_y|-D_y)C_F .
\end{equation}
Here \(D_yC_F=F'C_FD_z\); in particular the upper-right block is
\(C_F^\sharp|D_y|C_F\).
By \eqref{st-eq:congruence-D} and Remark~\ref{st-rem:H-not-small},
\begin{equation}\label{st-eq:H-full}
 \mathcal H_e=(|D_z|-D_z)+R_F .
\end{equation}
The smoothing kernel $R_F=C_F^\sharp|D_y|C_F-|D_z|$ is estimated in
Lemma~\ref{st-lem:kernel}. The order-one term $|D_z|-D_z$ retains
the flat contribution at every negative spatial index.
Both changes of spatial variable are static in the moving frame of
a fixed Stokes wave, so the conjugacy gives \eqref{st-eq:Atilde}
without a time derivative of the coordinate maps.

The coordinate and operator correspondence is
\begin{center}
\begin{tabularx}{\textwidth}{@{}l>{\raggedright\arraybackslash}X>{\raggedright\arraybackslash}X@{}}
\toprule
Coordinate & Linear operator & Change of variables \\
\midrule
$x=x_{\rm lab}-c_s t$ & $L_e^{\rm ph}$ & moving frame \\
$y$, with $x=X(y)$ & $L_e^{\rm c}$, or $A_e$ in complex notation
  & good unknown and conformal density \\
$z=F(y)$ & $\widetilde A_e$ & density with $F'=v/b$ \\
\bottomrule
\end{tabularx}
\end{center}

The table \eqref{st-eq:Y3} implies
\begin{equation}\label{st-eq:kap3}
 v=1-\tfrac12e^2+O(e^4),\qquad
 a-b=\tfrac12e^2+e^3C_1+O(e^4),\qquad
 \kappa(z)=1+2e^3\cos z+O(e^4).
\end{equation}
Indeed $1/b=1+2eC_1+O(e^2)$ gives
\[
 \frac ab=1+\frac{a-b}{b}=1+\tfrac12e^2+2e^3C_1+O(e^4),
 \qquad \frac{va}{b}=1+2e^3C_1+O(e^4).
\]
Since $F^{-1}=\Id+O(e)$, the cubic coefficient in the $z$ coordinate
is $2\cos z$.
Let
\[
 M(e)=\av{\kappa},\qquad W(e,z)=\kappa-M .
\]
We keep \(v\) and \(M\) exact. Their evenness follows from the
\(\pi\)-translation of the exact wave and the zero-mean primitive in
\eqref{st-eq:vF}. In particular \(M=1+O(e^4)\), and
\(W=2e^3\cos z+O(e^4)\).

\begin{remark}\label{st-rem:means}
The quantity $v$ is the harmonic mean of $b$ in the $y$ coordinate.
The quantity $M$ is the mean of $\kappa$ in the $z$ coordinate; changing
variables gives the useful identity
\[
 M=v^2\left\langle\frac{a}{b^2}\right\rangle_y.
\]
Thus the two coefficients of the spatially constant comparison law are
$v$ and $\sqrt M$. They are defined by this specific static reduction;
$v$ is distinct from the physical wave speed $c_s$.
\end{remark}

\subsubsection*{The physical momentum}
Using \(\Psi=c_s\mathcal H Y\) modulo constants gives
\begin{equation}\label{st-eq:action}
 I=-\av{\Psi Y'}=c_s\av{Y\K Y}
       =\tfrac12e^2+O(e^4).
\end{equation}
This is the physical momentum defined in \eqref{pc-momentum}.
It is even and analytic in \(e\). As a function of \(s=e^2\) its
derivative at zero is \(1/2\). The ordinary analytic inverse theorem
therefore gives \(e^2=2I+O(I^2)\).

\section{Sharp spectral comparison and sixth-order detuning}\label{st-sec:spectrum}
\subsection{Full operator families and the comparison order}
The part of \(\mathcal H_e\) at negative Fourier indices in
\eqref{st-eq:H-full} includes an order-one operator even at $e=0$.
In the first conformal coordinate $y$, let
$A_{\rightarrow}^{\rm c}=
\bigl(\begin{smallmatrix} D_y\circ b&D_y\\a&bD_y\end{smallmatrix}\bigr)$,
obtained from \eqref{st-eq:A} by replacing $|D_y|$ by $D_y$.
The density map \eqref{st-eq:density} and the identity
$C_F^\sharp D_y C_F=D_z$ transport it to
\[
 A_{\rightarrow}(e)=
       \begin{pmatrix}v(e)D&D\\\kappa(e)&v(e)D\end{pmatrix},
       \qquad \mathcal D(A_{\rightarrow})=H^1(\T;\mathbb C^2).
\]
Here and throughout this section, $D=D_z$. The comparison operator
acts on all spatial indices. Its flat blocks are
\[
 A_{\rightarrow}(0;n)=\begin{pmatrix}n&n\\1&n\end{pmatrix}.
\]
For \(n<0\) they have eigenvalues \(n\pm i\sqrt{|n|}\), whereas the
true flat blocks have \(n\pm\sqrt{|n|}\).
For \(n>0\) the two blocks agree.

Fix a nonsquare positive \(j\geq2\) and one sign \(\tau\).
The positive number \(\beta_0=j-\tau\sqrt j\) is a simple isolated
eigenvalue for both flat families. Their flat resolvents have
\(L^2\to H^1\) bounds on a common sufficiently small circle about it.
The coefficients \(v,\kappa\) are analytic in \(e\) on a fixed weaker
positive spatial strip. This follows from the exact Stokes branch,
scalar analytic operations, and the analytic inverse of the
near-identity map \(F\); it does not require analyticity of \(C_F\)
in \(\mathcal B(H^s,H^s)\).
Thus \(A_{\rightarrow}\) has an ordinary common-domain analytic
Riesz branch \(\beta_{\rightarrow}(e)\).
The same \(\pi\)-translation gives
\(A_{\rightarrow}(-e)=T_\pi^{-1}A_{\rightarrow}(e)T_\pi\);
hence this simple branch is even. For real $e$, reflection followed
by conjugation commutes with $A_{\rightarrow}$ since $\kappa$ is
real-valued and even in $z$. Simplicity then makes the selected eigenvalue real.
We set \(\lambda_{\rightarrow}=-\beta_{\rightarrow}\).
\begin{lemma}
\label{st-lem:transported-type-a}
After fixing a positive spatial strip and decreasing the amplitude radius,
the true transported family has the realization
\[
 \widetilde A_e=
 \begin{pmatrix}v(e)D&|D|+R_{F_e}\\
                 \kappa(e)&v(e)D\end{pmatrix},
 \qquad \mathcal D(\widetilde A_e)=H^1(\T;\mathbb C^2),
\]
and is analytic as a map into
\(\mathcal B(H^1,L^2)\). It has a rank-one analytic Riesz branch
at every nonsquare positive flat index under consideration.
For real amplitudes this branch equals the branch of the full physical
operator \(A_e\).
\end{lemma}
\begin{proof}
The Stokes branch, the scalar operations defining \(v\), and the inverse
circle map defining \(\kappa\) are holomorphic on a fixed weaker positive
strip once the amplitude disc has been reduced. The Hardy-kernel
construction in Section~\ref{st-sec:conformal} is likewise defined for
complex analytic \(F_e-\mathrm{Id}\) in an ordinary open ball.
The matrix bounds there give holomorphy of \(e\mapsto R_{F_e}\)
in \(\mathcal B(L^2,L^2)\) and on every fixed finite list of
Sobolev levels by the same weighted matrix convergence.
The displayed realization follows from
\(D+\mathcal H_e=|D|+R_{F_e}\).

It follows that \(\widetilde A_e-\widetilde A_0\) is small in
\(\mathcal B(H^1,L^2)\). The flat graph estimate and the resolvent
Neumann argument on a fixed small circle about the selected simple flat
eigenvalue give a common closed domain and a rank-one analytic Riesz
projection. For real \(e\), the genuine density map and its inverse
are bounded on both \(L^2\) and \(H^1\), and intertwine
\(A_e\) and \(\widetilde A_e\). Their selected branches are therefore
equal for real \(e\), hence equal as analytic germs.

The common-domain family \(\widetilde A_e\) thus supplies the
Taylor recursions directly.
\end{proof}

The true branch below is the common-domain branch of
Lemma~\ref{st-lem:transported-type-a}; for real amplitudes it is the
physical branch transported by \eqref{st-eq:density}.

\begin{lemma}
\label{st-lem:quadratic-transport}
On a fixed complex amplitude disc, $\|R_{F_e}\|_{\mathcal B(L^2,L^2)}
\le C|e|^2$, with $C$ independent of the selected spectral branch.
\end{lemma}
\begin{proof}
By \eqref{st-eq:Y3} and \eqref{st-eq:vF},
$F_e(y)=y+2e\sin y+O(e^2)$. Differentiate the congruence
$C_{F_e}^\sharp |D| C_{F_e}-|D|$ on Fourier polynomials at the identity.
The two composition terms together give
\[
 \left.\partial_e(R_{F_e})_{mn}\right|_{e=0}
 =i\bigl(n|m|-m|n|\bigr)\widehat{(2\sin y)}(m-n)=0.
\]
The prefactor vanishes for equal signs or a zero index. Opposite
nonzero signs require $|m-n|\ge2$, beyond the first harmonic of the
Stokes deformation. By the operator-valued holomorphy in
Lemma~\ref{st-lem:transported-type-a}, the bounded derivative is zero.
The value at zero is also zero, which gives the estimate.
\end{proof}

\begin{lemma}\label{st-lem:protection}
With \(j,\tau\) fixed as above, the true and comparison branches satisfy
\[
 \beta_{j,\tau}(e)-\beta_{\rightarrow,j,\tau}(e)
           =O_j(e^{2j+2}).
\]
The two branches here are those of the complete operators on $H^1$.
\end{lemma}
For a negative middle index $k$, a Fourier shift from $j$ to $k$
requires at least $j+|k|$ powers of $e$. The eigenvalue correction
requires a return to $j$, so the first possible order is $2j+2$.
The first possible contribution therefore has order six for $j=2$
and order at least eight for $j\ge3$. The proof uses a homotopy on the
full space and controls the two legs through right and left eigenvector
recursions. The constant potential is an absorbing direction in the
zero block.
\begin{proof}
Fix $j,\tau$, write $w=\tau\sqrt j$ and $\beta_0=j-w>0$, and
introduce the auxiliary family
\begin{equation}\label{st-eq:homotopy}
 A_\alpha(e)=
 \begin{pmatrix}v(e)D&D+\alpha\mathcal H_e\\
                 \kappa(e)&v(e)D\end{pmatrix},
 \qquad 0\leq\alpha\leq1,
 \qquad \mathcal D(A_\alpha)=H^1(\T;\mathbb C^2).
\end{equation}
Only the endpoints are the comparison and transported physical
operators. By \eqref{st-eq:H-full}, this is a common-domain analytic
family in $e$. Its distributional maximal $L^2$ domain is also
$H^1$: the second row gives $vDp\in L^2$, and then the order-one
upper-right entry applied to $p$ is in $L^2$, giving $vD\xi\in L^2$
from the first row. Its flat block at index $n$ is
\[
 A_{\alpha,0}(n)=\begin{pmatrix}n&a_\alpha(n)\\1&n\end{pmatrix},
 \qquad a_\alpha(n)=n+\alpha(|n|-n).
\]
For $n=-N<0$ its eigenvalues are
$-N\pm\sqrt{N(2\alpha-1)}$ and have real part at most zero.
The zero block is $\bigl(\begin{smallmatrix}0&0\\1&0\end{smallmatrix}\bigr)$.
The positive blocks do not depend on $\alpha$. Hence one circle about
$\beta_0$ isolates the same simple positive eigenvalue for every
$\alpha\in[0,1]$.

On that circle the flat resolvent at large $|n|$ is $O(|n|^{-1})$,
uniformly in $\alpha$: its determinant is $n^2+O(|n|)$ and its
numerator is $O(|n|)$. The finitely many remaining inverses are
continuous on a compact set without poles. Thus the flat resolvents
are uniformly bounded from $L^2$ to $H^1$.
The perturbation is $O(e)$ from $H^1$ to $L^2$, uniformly in
$\alpha$, so a common small amplitude disc gives a rank-one analytic
Riesz branch $\beta_\alpha(e)$ and a common closed domain.
For spectral parameter $\zeta$ on this contour, differentiation in
$\alpha$ gives
$\partial_\alpha(A_\alpha-\zeta)^{-1}
=-(A_\alpha-\zeta)^{-1}(\partial_\alpha A_\alpha)(A_\alpha-\zeta)^{-1}$.
The two resolvents map $L^2$ to $H^1$, and the middle factor maps
$H^1$ to $L^2$. The contour formula and the rank-one normalization
below give $\partial_\alpha\beta_\alpha$ holomorphic in $e$ and
uniformly bounded on a common smaller disc, including the one-sided
derivatives at the endpoints.
At index $j$ choose the flat right vector and left row
\[
 r_0=\binom{-w}{1}e^{ijz},\qquad
 \ell_0u=-\widehat u_1(j)/(2w)+\widehat u_2(j)/2,
 \qquad \ell_0r_0=1.
\]
For the comparison below fix $R=2j+2$ and $m=2R+4$ first.
On each integer level $H^s$, $|s|\leq m$, with domain $H^{s+1}$,
the same Fourier inverse and Neumann argument apply, after a single
further reduction of the radius. The flat Fourier inverse and each
factor of the Neumann series agree on common data. The series converge
in their respective Sobolev norms and hence in distributions, so the
resolvents and the Riesz projection $\mathcal P_\alpha(e)$ agree on
the finite scale. With
$d_\alpha=\ell_0\mathcal P_\alpha r_0$, the normalizations used below are
\[
 u_\alpha=\mathcal P_\alpha r_0/d_\alpha,\qquad
 \ell_\alpha=\ell_0\mathcal P_\alpha/d_\alpha .
\]
These are legitimate since $d_\alpha(0)=1$. They separately fix
$\ell_0u_\alpha=1$ and $\ell_\alpha r_0=1$; in general
$\ell_\alpha u_\alpha=1/d_\alpha$, which is retained in the
eigenvalue derivative below.
The positive levels give regularity of $u_\alpha$; the negative
levels make $\ell_\alpha$ bounded on $H^{-m}$, which gives the
corresponding positive regularity of its Fourier coefficient row.
Only this fixed finite reserve is used.

For the explicit coefficient recursion put $\Pi_0=r_0\ell_0$,
$Q_0=\Id-\Pi_0$, and let $S_\alpha$ be the
inverse $Q_0L^2\to Q_0H^1$ of $A_{\alpha,0}-\beta_0$ on the
complement, extended by zero on $\Pi_0L^2$. It is diagonal in the
spatial index. For $n\ne j$ its complete block is
\begin{equation}\label{st-eq:homotopy-resolvent}
 S_\alpha(n)=
 \frac1{(n-\beta_0)^2-a_\alpha(n)}
 \begin{pmatrix}n-\beta_0&-a_\alpha(n)\\-1&n-\beta_0\end{pmatrix}.
\end{equation}
In particular the two endpoint inverses agree for $n\geq0$ and
generally differ for $n<0$. At $j$, the reduced block retains the
other eigenvector of the same two-component Fourier block.

Write $A_\alpha=A_{\alpha,0}+\sum_{r\geq1}e^rV_{\alpha,r}$,
$\beta_\alpha=\beta_0+\sum_{n\geq1}e^n\beta_{\alpha,n}$, and
$u_\alpha=r_0+\sum_{n\geq1}e^nu_{\alpha,n}$, with
$\ell_0u_{\alpha,n}=0$. Comparing coefficients gives
\begin{align}
 \beta_{\alpha,n}
  &=\ell_0\sum_{r=1}^nV_{\alpha,r}u_{\alpha,n-r},
       \label{st-eq:right-recursion}\\
 u_{\alpha,n}
  &=-S_\alpha\left\{
       \sum_{r=1}^nV_{\alpha,r}u_{\alpha,n-r}
       -\sum_{r=1}^{n-1}\beta_{\alpha,r}u_{\alpha,n-r}\right\}.
       \nonumber
\end{align}
The missing term $\beta_{\alpha,n}r_0$ is annihilated by $S_\alpha$.
For the analytic left eigenrow
$\ell_\alpha=\ell_0+\sum_{n\geq1}e^n\ell_{\alpha,n}$ normalized by
$\ell_\alpha r_0=1$, the corresponding identity is
\[
 \ell_{\alpha,n}=-\left\{
    \sum_{r=1}^n\ell_{\alpha,n-r}V_{\alpha,r}
    -\sum_{r=1}^{n-1}\beta_{\alpha,r}\ell_{\alpha,n-r}
                          \right\}S_\alpha.
\]
Left rows are complex-linear Fourier functionals, not Hilbert
adjoints taken at a complex amplitude. All products in these
recursions have their indicated $H^1$ and $L^2$ meanings.

Each $V_{\alpha,r}$ changes a spatial label by at most $r$.
For its upper-right entry this follows by expanding
$\mathcal H_e=C_{F_e}^\sharp(2|D|P_-)C_{F_e}$ and applying
Lemma~\ref{st-lem:filtration} to both outer factors. The adjoint
matrix entries satisfy $(C_F^\sharp)_{mn}=(C_F)_{-n,-m}$ by
Appendix~\ref{st-sec:conformal}, so their Taylor orders obey the
same bandwidth bound. Here $P_-$ projects onto the negative
Fourier indices. Since $S_\alpha$ does not change
labels, induction in the right and left recursions proves
\[
 \widehat u_{\alpha,n}(k)=0,\qquad
 \ell_{\alpha,n}\big|_{\mathbb C^2e^{ikz}}=0
       \quad\hbox{if }|k-j|>n.
\]
The scalar eigenvalue feedback terms are included in this induction.

For real small $e$, differentiation of the simple eigenvalue in
$\alpha$ gives
\begin{equation}\label{st-eq:homotopy-derivative}
 \partial_\alpha\beta_\alpha(e)=
 \frac{\ell_{\alpha}^{[1]}C_{F_e}^\sharp
              (2|D|P_-)C_{F_e}u_{\alpha}^{[2]}}
      {\ell_\alpha u_\alpha}.
\end{equation}
Here bracketed superscripts $[1],[2]$ indicate components, whereas
the subscripts in the recursions indicate Taylor order. The denominator
equals one at $e=0$ and stays nonzero. For a negative Fourier index $k$ in the intermediate coordinate $y$,
put
\[
 a_k(e)=\widehat{C_{F_e}u_{\alpha}^{[2]}}(k),\qquad
 b_k(e)=\ell_{\alpha}^{[1]}C_{F_e}^\sharp e^{iky}.
\]
The filtration and the triangle inequality show that each factor
has zero Taylor coefficients below degree $|k-j|=j+|k|$.
The numerator in \eqref{st-eq:homotopy-derivative} is therefore
$2\sum_{k<0}|k|b_k(e)a_k(e)$. For real amplitudes this sum is
absolutely convergent: $C_{F_e}u_\alpha^{[2]}\in H^1$, and the
row $\ell_\alpha^{[1]}C_{F_e}^\sharp$ is bounded on $L^2$, so
Cauchy--Schwarz applies to $\{|k|a_k\}$ and $\{b_k\}$.
The coefficient filtration shows that the sum has no coefficient below
degree $2j+2$. A coefficient of degree $R$ can involve only
$2(j+|k|)\leq R$, so only finitely many negative middle indices
contribute.

For the regularity of this extraction, the chain rule for real $e$
writes $\partial_e^r C_{F_e}f$ on smooth functions as a sum of smooth
coefficients times $C_{F_e}\partial_z^k f$, $0\le k\le r$.
Change of variables and the spatial chain rule give
\[
 \|\partial_e^r C_{F_e}f\|_{H^s}\le C_{r,s}\|f\|_{H^{s+r}},
 \qquad s=0,1,\quad 0\le r\le R,
\]
uniformly on a compact small real amplitude interval with $F'_e>1/2$.
Density and the uniform bound extend the derivatives and their strong
continuity to fixed data at the stated levels. In the left row use
the bilinear pairing $\ell_\alpha^{[1]}f=\langle g f\rangle$,
where $g\in H^m$ by boundedness of the row on $H^{-m}$. On the real
slice $\ell_\alpha^{[1]}C_{F_e}^\sharp f=\langle(C_{F_e}g)f\rangle$.
The high-level analyticity of the vector and the row and the chain-rule
bounds give $C^R$ factors in $H^1$ and $L^2$, respectively.
Differentiate their continuous bilinear pairing across
$|D|P_-:H^1\to L^2$. Each resulting Taylor coefficient has finite
Fourier support by the eigenvector and composition bandwidth bounds,
so the coefficient sum is finite.

Zero modes have not been removed from the argument. The first output
component has zero mean, and a constant second component is killed
by every $A_\alpha(e)$. At index zero,
\[
 S_\alpha(0)=
 \begin{pmatrix}-\beta_0^{-1}&0\\-\beta_0^{-2}&-\beta_0^{-1}\end{pmatrix}.
\]
Thus a constant second component stays in that absorbing sector;
no first-component mean can be generated from the selected mode.
Also $|D|P_-$ has no zero-mode contribution. No shorter return
through zero is possible.

All estimates are uniform on the compact $\alpha$ interval.
Analyticity on a common amplitude disc now gives
$\partial_\alpha\beta_\alpha=O_j(e^{2j+2})$ uniformly there.
Integrating from $\alpha=0$ to $1$ gives the claimed comparison bound.
\end{proof}

\begin{remark}
On positive inputs the entries of \(\mathcal H_e\) enjoy both the
analytic smoothing estimate and a Taylor-order restriction.
The all-order Fourier filtration of the actual Stokes branch
supplies the factor \(e^{2j+2}\).
\end{remark}

\subsection{Sharpness of the comparison order}\label{st-sec:sharp-protection}
The leading correction uses only the highest harmonics of the Stokes
profile and the flat inverse at index $-1$.

\begin{lemma}[Highest harmonics]\label{st-lem:highest}
For a Taylor--Fourier series $f$ with bandwidth at most $r$ at
degree $r$, write $f_r=[e^r]f$ and
$f^\natural(\varrho)=\sum_{r\ge0}\widehat f_r(-r)\varrho^r$ in a
formal variable $\varrho$. This series records the most negative
allowed Fourier index at each Taylor degree.
For the Stokes coefficients of Section~\ref{st-sec:stokes},
$T=2Y^\natural$ satisfies $T=\varrho\exp(T)$, the normalized coefficient satisfies
$\kappa^\natural=1$, and
$[e^{j+1}](C_F)_{-1,j}=-j/(j+1)!$ for $j\ge2$, where
$(C_F)_{-1,j}$ is the Fourier coefficient at $-1$ of $C_F(e^{ijz})$.
In particular, writing $\kappa_r=[e^r]\kappa$, one has
$\kappa_1=0$ and $\operatorname{supp}\widehat\kappa_r\subset
\{|k|\le r-2\}$ for $r\ge2$.
\end{lemma}
\begin{proof}
The identities in this proof are coefficient identities in the
formal Taylor--Fourier algebra of Lemma~\ref{st-lem:filtration}.
At each degree the products and compositions use only finitely many
terms. The analytic wave and the finite-derivative composition
estimates identify the extracted coefficients with those of the
actual maps.
If a product contributes at degree $r$ and index
$-r$, each factor is at its own most negative index. Thus this
operation preserves products and inverses of units. It sends $|D_y|$
to $\varrho\partial_\varrho$, written $\mathscr D$, and $\partial_y$ to
$-i\mathscr D$. In particular, any scalar function of $e$ contributes
only its value at zero.

Put $T=2Y^\natural$. The Babenko identity \eqref{st-eq:Babenko}
gives $(1-T)\mathscr DT=T$ with $T=\varrho+O(\varrho^2)$. Equivalently
$\mathscr D(T\exp(-T))=T\exp(-T)$, so $T=\varrho\exp(T)$.
The conformal derivative is supported on nonpositive Fourier indices,
with $(Z')^\natural=1+\mathscr DT$. Reflecting the other factor in
$J=|Z'|^2$ gives $J^\natural=1+\mathscr DT=(1-T)^{-1}$.
Consequently the coefficient identities in Section~\ref{st-sec:physical}
give
\[
 b^\natural=1-T,\qquad B^\natural=\tfrac i2T,
 \qquad (bB_y)^\natural=\tfrac12T,\qquad
 a^\natural=\frac{1-T/2}{1+\mathscr DT/2}=1-T.
\]
Here $B=-bY'$. Since $F'=v/b$ and the primitive has zero mean,
$(F-y)^\natural=iT$. For composition with $F$, the maximal harmonic
of $g\circ F$ is the series $g^\natural(\varrho\exp(-i(F-y)^\natural))$.
Indeed an input at order $r$ can reach a maximal output harmonic only
if its own index is $-r$ and the subsequent shifts are maximal.
Now $(\kappa\circ F)^\natural=(va/b)^\natural=1$ and
$\varrho\exp(-i(F-y)^\natural)=\varrho\exp(T)=T$.
Since $T$ has linear coefficient $1$, it follows that
$\kappa^\natural=1$. The evenness of $a,b$ and the oddness of $F$
make $\kappa$ even in $z$, ruling out the positive maximal harmonics.
The relation
$\kappa(-e,z)=\kappa(e,z+\pi)$ makes its harmonic indices at degree
$r$ congruent to $r$ modulo two. This gives the stated bandwidth.
By Lagrange inversion, with $\upsilon$ the coefficient variable on the right,
\[
 [e^{j+1}](C_F)_{-1,j}=[\varrho^{j+1}]\exp(-jT(\varrho))
 =-\frac{j}{j+1}[\upsilon^j]\exp(\upsilon)=\gamma_j .
\]
\end{proof}

\begin{proof}[Proof of Theorem~\ref{st-thm:sharp-protection}]
Use the rows, vectors and coefficients $a_k,b_k$ of the proof of
Lemma~\ref{st-lem:protection}. At degree $2j+2$, only the middle
index $k=-1$ can contribute to \eqref{st-eq:homotopy-derivative}.
For a path $j=k_0,\ldots,k_p=-1$ contributing to a degree-$j+1$ leg,
the orders $r_\nu$ of its factors satisfy
\[
 j+1=\sum_{\nu=1}^p r_\nu\ge
     \sum_{\nu=1}^p|k_\nu-k_{\nu-1}|\ge j+1.
\]
Equality forces monotone motion, with each shift made at its
lowest possible order. In the left recursion we read the path from
the selected row to its successive column labels, so the physical
matrix entries occur with inputs and outputs reversed. The order
bounds below are symmetric under this reversal.
Positive-order terms that leave the label fixed, including scalar
eigenvalue insertions, cannot occur.
The maximal-harmonic vanishing just proved excludes a tight shift
through the lower-left coefficient $\kappa$. The upper-right
perturbation between positive indices $m,n$ has degree at least
$m+n+2$. In \eqref{st-eq:hardy} the summand at $k=-N$ requires
degrees at least $m+N$ and $n+N$ in the two composition factors.
The zero row and column of $R_F$
vanish, and a constant potential has no outgoing term. Thus there
is no tight intermediate step between distinct
nonnegative indices. The only possible tight move within either
eigenvector recursion is a jump directly from $j$ to $-1$ along
its path, using the upper-right entry with the corresponding
input-output order. Otherwise the entire shift occurs in the
outer composition factor.

At degree $j+1$ that entry is $2\alpha\gamma_j$ in both directions.
In \eqref{st-eq:hardy}, the index $-1$ already costs $j+1$ on the
positive leg, and the other composition is the identity at degree
zero. The real transport map is odd, so the composition coefficients
are real, and $R_F$ is symmetric for real $e$. This gives the reverse
entry as well. In the final composition factors of $a_{-1}$ and
$b_{-1}$, only the direct composition from $j$ and the eigenvector
coefficient at $-1$ can contribute at degree $j+1$.

Write $D_\alpha=(j+1-w)^2+1-2\alpha$. The $(2,1)$ entry of the
shifted reduced inverse at $-1$ is $-1/D_\alpha$ by
\eqref{st-eq:homotopy-resolvent}. In the right recursion the
contribution is $-S_\alpha(-1)(2\alpha\gamma_j E_{12})r_0(j)$,
where $E_{12}$ has its only unit entry at $(1,2)$.
In the left recursion it is
$-\ell_0(j)(2\alpha\gamma_j E_{12})S_\alpha(-1)$.
Here $r_0(j)$ and $\ell_0(j)$ are the flat column and row at
Fourier index $j$. Taking the second component of the column and
the first entry of the row, then adding the direct composition
terms, gives
\[
 [e^{j+1}]a_{-1}=\gamma_j\Bigl(1+\frac{2\alpha}{D_\alpha}\Bigr),
 \qquad
 [e^{j+1}]b_{-1}=-\frac{\gamma_j}{2w}
                 \Bigl(1+\frac{2\alpha}{D_\alpha}\Bigr).
\]
The denominator in \eqref{st-eq:homotopy-derivative} has constant
coefficient $1$. The coefficient of its numerator at degree $2j+2$
is thus $-\gamma_j^2(D_0/D_\alpha)^2/w$. Integrating over $\alpha$
from $0$ to $1$ uses $D_\alpha=D_0-2\alpha$ and
$\int_0^1(D_0/D_\alpha)^2\,d\alpha=D_0/D_1$.
This ratio accounts for the contribution propagated through the
flat block at $-1$. The direct composition by itself gives the
term without the factor $D_0/D_1$.
The common-domain analytic construction and the finite derivative
reserve in Lemma~\ref{st-lem:protection} justify the coefficient
calculation. The half-period symmetry removes odd orders. Changing
from $\beta$ to $\lambda$ gives \eqref{st-eq:sharp-protection}.
Since $j+1-w>1$ for $j\ge2$, the coefficient has sign $\tau$.
The denominator in the second summand grows like $j^2$
for either sign, giving its stated large-$j$ asymptotic.
Substitution at $j=2$ gives \eqref{st-eq:boundary}.
The momentum version uses
$(e^2)^{j+1}=(2I+O(I^2))^{j+1}$ at each fixed $j$.
\end{proof}

\subsection{Two cubic insertions and sixth-order detuning}\label{st-sec:return}
Use the comparison family only to compute the protected coefficients.
Keep the mean part exact:
\[
 \overline A(e)=\begin{pmatrix}vD&D\\M&vD\end{pmatrix},\qquad
 A_{\rightarrow}-\overline A=\begin{pmatrix}0&0\\W&0\end{pmatrix}.
\]
Its branch with sign $\tau$ at positive Fourier index \(j\) is
\[
 \overline\beta_{j,\tau}=jv-\tau\sqrt{jM}.
\]
Every first-order insertion of \(W\) into this spatially diagonal
problem has zero scalar projection, because \(\av W=0\) exactly.
Since \(W=e^3(2\cos z)+O(e^4)\), the only nonconstant contribution
of degree six is the double insertion of \(W_3=2\cos z\).
The mean coefficients \(v-1=O(e^2)\), \(M-1=O(e^4)\) can be set
to zero in this degree-six double insertion: their effects start
at degree eight.

Put \(w=\tau\sqrt j\) and choose
\[
 r=\binom{-w}{1},\qquad
 l=\begin{pmatrix}-1/(2w)&1/2\end{pmatrix},\qquad lr=1.
\]
The offsets of $W_3$ are $h=\pm1$, each with matrix
$V_h=\bigl(\begin{smallmatrix}0&0\\1&0\end{smallmatrix}\bigr)$.
For $j\ge2$ both intermediate indices $j+h$ are positive. The
contribution at offset $h$ to $\beta$ is
$-lV_{-h}(A_{\rightarrow}(0;j+h)-\beta_0)^{-1}V_h r$.
The full two-component denominator is
\[
 \det\begin{pmatrix}h+w&j+h\\1&h+w\end{pmatrix}
       =h(h+2w-1).
\]
The full inverse includes both eigenvectors of the intermediate Fourier block. With the sign
from the reduced resolvent, the eigenvalue correction is
\begin{align}
 B_6(w)
 &=-\frac w2\sum_{h=\pm1}\frac{j+h}{h(h+2w-1)} \nonumber\\
 &=-\frac w2\left\{\frac{w^2+1}{2w}
                        +\frac{w^2-1}{2(1-w)}\right\}
   =\frac{w-1}{4}. \label{st-eq:return}
\end{align}
For the nonsquare indices at issue these denominators are nonzero.
\begin{proof}[Proof of Proposition~\ref{st-prop:individual}]
The sign change from $\beta$ to $\lambda$ gives the comparison correction
$(1-w)e^6/4$. Parity removes odd orders. For $j\ge3$,
Theorem~\ref{st-thm:sharp-protection} transfers the expansion to
the physical branch. The comparison expansion itself also holds at $j=2$.

For the variance form, expand the square root in $d_s$ with
$\kappa=M+W$ and $\langle W\rangle=0$:
\[
 d_s=\sqrt M-\frac{\langle W^2\rangle}{8M^{3/2}}+O(e^9)
     =\sqrt M-\frac14e^6+O(e^8),\qquad
 m_s=\frac14e^6+O(e^8).
\]
Here $\langle W^2\rangle=2e^6+O(e^7)$; half-period translation makes
$d_s$ and $m_s$ even, removing the seventh-order terms. Substitution
gives \eqref{st-eq:variance-expansion}.
\end{proof}

\begin{proof}[Proof of Theorem~\ref{st-thm:detuning}]
For \(j=4q,q,q\) with signs \(+,-,+\), the mean transport coefficient
in the combination is
\[
 -4q+q+3q=0,
\]
and the complete mean dispersion coefficient is
\[
 2\sqrt q+\sqrt q-3\sqrt q=0.
\]
These cancellations are identities in the exact functions \(v(e)\)
and \(\sqrt{M(e)}\), not cancellations in a finite truncation.
The remaining comparison coefficient at order six is $-3/4$.
For $q\ge3$ this is also the physical coefficient by
Lemma~\ref{st-lem:protection}. For $q=2$, the same comparison
calculation uses the positive blocks $1,3$ for each branch at $j=2$.
The two shifted block determinants are $2(1-w)$ and $2w$.
They are nonzero, so the additional physical coefficients come from
\eqref{st-eq:boundary}, with the branch at $j=8$ still protected.
The amplitude coefficient is
\[
 -\frac34-\frac{1-2\sqrt2}{21}-3\frac{1+2\sqrt2}{21}
 =-\frac{79+16\sqrt2}{84}.
\]
In each case $e^6=8I^3+O(I^4)$ by \eqref{st-eq:action}.
Even analyticity in \(e\) gives analyticity in \(I\), hence the
derivative and sign estimates in the theorem.
\end{proof}

\section{Eighth-order detuning}\label{st-sec:eighth}

The eighth-order relation has two contributions. The differential
comparison acquires a mode-dependent correction even when the means
are kept exact. The branch at $j=3$ also detects the negative
spatial indices at that order.

\subsection{The fourth-order nonconstant coefficient}
Write $K=\kappa\circ F$, $C_n=\cos(ny)$ and $S_n=\sin(ny)$.
A subscript on a spatial coefficient denotes its Taylor degree in $e$.
To compute $Y_4$, use $Y_1,Y_2,Y_3$ from \eqref{st-eq:Y3}
and $\sigma_2=1$ in \eqref{st-eq:recursion}. Product-to-sum identities
give
\[
 \begin{aligned}
 \sum_{a+b=4}Y_a|D|Y_b
   &=6C_1C_3+2C_2^2-C_2=1+2C_2+4C_4,\\
 \tfrac12|D|\sum_{a+b=4}Y_aY_b
   &=\tfrac12|D|(3C_1C_3+\tfrac14-C_2+C_2^2)
      =\tfrac12C_2+4C_4 .
 \end{aligned}
\]
The odd speed coefficients vanish. Subtracting $|D|Y_2=2C_2$
and dividing by the symbol of $|D|-1$ at indices $0,2,4$ yields
$Y_4=-1+\tfrac12C_2+\tfrac83C_4$.

Since $X'=1+|D|Y$,
\[
 J_4=2|D|Y_4+2\big((|D|Y_1)(|D|Y_3)+Y_1'Y_3'\big)
                    +(|D|Y_2)^2+(Y_2')^2.
\]
The last two terms sum to $4$ and the cross term is $9C_2$, giving
$J_4=4+11C_2+\tfrac{64}{3}C_4$. For $B=-bY'$, the first three coefficients
of $b$ in \eqref{st-eq:Y3} give
\[
 B_4=S_2+\tfrac{32}{3}S_4-9C_1S_3
           +2(\tfrac32-2C_2)S_2+(C_1-3C_3)S_1.
\]
Reducing the products gives $B_4=\tfrac32S_2+\tfrac83S_4$.

The identities $b=c_s/J$ and $a=(1-bB_y)/X'$ give
$KX'=(v/c_s)J-vB_y$. Since $v/c_s=\langle J\rangle_y^{-1}$ by
\eqref{st-eq:vF}, the means of $J_2$ and $J_4$ give
$v/c_s=1-e^2-3e^4+O(e^6)$. Using $v=1-\tfrac12e^2+O(e^4)$,
$K=1+2e^3 C_1+O(e^4)$ and the lower coefficients from \eqref{st-eq:Y3},
we obtain
\[
 \begin{aligned}
 [e^4]((v/c_s)J)&=7C_2+\tfrac{64}{3}C_4,\\
 [e^4](vB_y)&=2C_2+\tfrac{32}{3}C_4,\\
 [e^4](KX')&=K_4+1+2C_2+\tfrac{32}{3}C_4 .
 \end{aligned}
\]
In the last line the cross term is $K_3X'_1=2C_1^2$.
Thus $K_4=-1+3C_2$. The inverse coordinate
$F^{-1}(z)=z-2e\sin z+O(e^2)$ adds $4\sin^2 z$ to the fourth-degree
coefficient when applied to $2e^3\cos y$. This gives
$\kappa_4(z)=1+\cos 2z$. Subtracting the mean yields
\[
 W(e,z)=2e^3\cos z+e^4\cos 2z+e^5 W_5(z)+O(e^6),
 \qquad W_5=[e^5]W.
\]

\subsection{The second return with exact means}
Fix $j\ge3$ nonsquare and $w=\tau\sqrt j$. The exact mean operator $\overline A$
in Section~\ref{st-sec:return} has eigenvalue
$\overline\beta=jv-w\sqrt M$ and right and left eigenvectors
at index $j$ given by
\[
 r=\binom{-w/\sqrt M}{1},\qquad
 l=\begin{pmatrix}-\sqrt M/(2w)&1/2\end{pmatrix}.
\]
Here $lr=1$. Let $P$ be the corresponding rank-one spectral projection,
$P_c=\Id-P$, and $V$ the perturbation with sole entry
$V_{21}=W$. These projections commute with the mean operator.
The reduced mean inverse is bounded from $P_c L^2$
to $P_c H^1$ for each fixed $j$ and small $e$. This follows
from the flat gap and the common-domain estimate in
Lemma~\ref{st-lem:protection}, or by adding $P$ to the shifted
operator and using the bounded inverse theorem there.
Let $\beta_\rightarrow$ denote the selected eigenvalue of
$A_\rightarrow$. Its displacement from $\overline\beta$ is
$O_j(e^3)$ by the rank-one Riesz formula since $V=O(e^3)$ on $L^2$.
In the exact complementary inverse, dropping the perturbation
$P_c\bigl(V-(\beta_\rightarrow-\overline\beta)\Id\bigr)P_c$
changes the inverse by $O_j(e^3)$ from $L^2$
to $H^1$. The two outer $V$ factors make the resulting error
$O_j(e^9)$. The direct term $lVr$ vanishes since $W$ has zero mean.

In the shifted block at $j+h$, $h\ne0$, the determinant is
$h(v^2h+2vw\sqrt M-M)$. Substituting its $(1,2)$ inverse entry and
using $\lambda=-\beta$ yields
\begin{equation}\label{st-eq:eighth-schur}
 \lambda_{\rightarrow,j,\tau}-(-jv+w\sqrt M)
 =\frac{w}{2\sqrt M}\sum_{h\ne0}
   \frac{(j+h)\widehat W(e,h)\widehat W(e,-h)}
        {h(v^2h+2vw\sqrt M-M)}+O_j(e^{10}) .
\end{equation}
The sum is the bounded reduced-resolvent pairing just obtained.
In particular the infinite sum is meaningful before extracting
coefficients. Both the left side and the sum are even analytic
functions of $e$ by half-period translation, which removes the
order-nine error. The zero Fourier block corresponds to $j+h=0$
and contributes zero.

Write $d_h=h(h+2w-1)$ for the flat value of that shifted determinant. The denominator is
$d_h-e^2h(h+w)+O(e^4)$ and $M^{-1/2}=1+O(e^4)$.
Put $a_5=\widehat W_5(1)=\widehat W_5(-1)$, with equality by
reflection. At order eight, only $h=\pm1,\pm2$ occur.
The $W_3W_5$ product contributes $2a_5(1-w)/4$ by
\eqref{st-eq:return}. Differences of the normalized gaps therefore require
only the $W_4W_4$ product and the order-two change in the denominators
of the $W_3W_3$ product. The pairings are
\begin{align*}
 \frac w2\sum_{h=\pm1}\frac{(j+h)h(h+w)}{d_h^2}
   &=\frac{w+1}{8w},\\
 \frac w8\sum_{h=\pm2}\frac{j+h}{d_h}
   &=-\frac{w(w-1)^2}{4(2w-3)(2w+1)} .
\end{align*}
For $j=3$ the term at $h=-2$ visits index $1$. At $e=0$ its shifted determinant is
$-2(2w-3)\ne0$ for $w=\pm\sqrt3$.
For example the second sum uses
$(w^2+2)/(2w+1)-(w^2-2)/(2w-3)$ after factoring out $w/16$.
Let
\[
 g(w)=\frac{w+1}{8w}-\frac{w(w-1)^2}{4(2w-3)(2w+1)} .
\]
Subtracting
the values at $w$ and $-w$ gives
\[
 g(w)-g(-w)=
 \frac{-8j^3+30j^2-34j+9}{4w(4j-1)(4j-9)},\qquad j=w^2.
\]

\begin{proof}[Proof of Corollary~\ref{st-cor:eighth-detuning}]
The four weights annihilate the exact mean terms and every affine
function of $w$. In particular they annihilate the sixth-order
return and the order-eight term containing $a_5$.
At $w=\sqrt3$ and $w=2\sqrt3$, the difference $g(w)-g(-w)$ equals
$-13/(44\sqrt3)$ and
$-3301/(4888\sqrt3)$, respectively. Twice the first minus the second is
$1513\sqrt3/53768$, the order-eight coefficient for the comparison.
Theorem~\ref{st-thm:sharp-protection} gives the coefficient
$(2+7w)/1056$ at $j=3$ and no term of order eight at $j=12$.
Thus the exact physical contribution from the negative indices is
$7\sqrt3/264$. Since the degree-six term cancels, changing variables
uses only $e^8=16I^4+O(I^5)$ and proves \eqref{st-eq:eighth}.
\end{proof}

\begin{proof}[Proof of Theorem~\ref{st-thm:gap-variation}]
For the comparison, the contribution of $\overline A$ to the
normalized gap is $\sqrt M$. The degree-six correction $-e^6/4$
is independent of $j$. The order-eight
term containing $a_5$ is also common. Dividing $g(w)-g(-w)$ by $2w$
gives $f(j)-1/16$. Both signs at every nonsquare $j\ge5$
agree with the physical branches through order eight by
Theorem~\ref{st-thm:sharp-protection}. At $j=3$ the same theorem
adds $(C_{3,+}-C_{3,-})/(2\sqrt3)=7/1056$ to the coefficient of
$e^8$ in the normalized gap. The common terms, including the one
with $a_5$, cancel in the difference. Now convert $e^8$ to $I^4$.
For $x\ge4$,
\[
 f'(x)=-\frac1{8x^2}+\frac9{16(4x-1)^2}+\frac3{16(4x-9)^2}
 \le\frac{-1/8+1/25+3/49}{x^2}<0 .
\]
The boundary comparison is $\mathfrak f(3)=7/352$ and
$\mathfrak f(3)-f(5)=219/33440>0$. Thus the coefficient
has the sign of $k-j$ for distinct nonsquare $j,k\ge3$.
\end{proof}

\begin{remark}[Comparison with the normal-form frequencies]
\label{st-rem:normal-form}
In the orthonormal Fourier convention of Feola--Giuliani, the momentum
in their complex variables is $-\sum_n n|u_n|^2$ and is preserved
by their Birkhoff map~\cite[(2.1), (3.8), (4.6)]{FG21}. Our $I$ is
minus that momentum divided by $2\pi$. On the reference torus with
tangential set $\{1\}$ and action $\mathcal J=\varepsilon^2\zeta_1$,
this gives $I=\mathcal J/(2\pi)$~\cite[(5.5)]{FG21}. In comparing actual
profiles one must use the full momentum, including corrections to
the reference torus. Our amplitude is fixed independently by the
first conformal cosine and \eqref{st-eq:action}.

In their doubled complex convention the final diagonal entries at
spatial index $n$ have signed frequencies $d_n^\infty$ and
$-d_{-n}^\infty$ (arXiv v3, (12.2), (12.5)). The conversion to our
frequencies would require identifying the exact Stokes background and
the corresponding spectral subspaces of the full physical operator.
For this conversion, suppose a conjugacy $T(\theta)$ from reduced to
physical variables satisfies the single-carrier translation covariance
$T(\theta-s)U_s=U_sT(\theta)$, where $U_s f(x)=f(x+s)$.
This is the symmetry in arXiv v3, (1.22) and Lemma~2.22.
Then $U_\theta T(\theta)U_{-\theta}=T(0)$. At $\theta=c_s t$
the transformation in the moving frame is constant, and translation
of the reduced mode supplies $e^{inc_st}$. The frequency dictionary
under these identifications is
\[
 \lambda_{n,+}=d_n^\infty-nc_s,\qquad
 \lambda_{n,-}=-d_{-n}^\infty-nc_s .
\]
The transport contribution $-nc_s$ cancels in $\Delta_q$.
In their formula (11.3) for the intermediate diagonal entries $d_n$,
$\mathtt m_0^{(\ge4)}$ is the coefficient of $\operatorname{sign}(n)$,
and $r_n$ is the mode-dependent remainder. Their combination is
\[
 d_{4q}+d_{-q}-3d_q=-3\mathtt m_0^{(\ge4)}+r_{4q}+r_{-q}-3r_q .
\]
The spatial decay estimate (11.4) leaves the fixed-index sixth-order
coefficients of $r_n$ undetermined. They include diagonal terms from
(11.79) and Section~11.4. For the correction $r_n^{(\infty)}$ from the intermediate to the final
diagonal, (12.3) of arXiv v3 bounds $\sup_n\langle n\rangle^{1/2}
|r_n^{(\infty)}|^{\gamma^3}$ by $C\varepsilon^{13}\gamma^{-3}$.
In their notation, $|f|^{\gamma^3}$ is the norm
$\sup|f|+\gamma^3\Lip(f)$ on the parameter set (Section~2.1).
It bounds the unweighted supremum directly. In the regime
(13.1)--(13.2), $\gamma=\varepsilon^{2+a}$ and $0<a<1/8$, so the
pointwise error is $O(\varepsilon^{7-3a})$. On a matched one-mode
family with $e\asymp\varepsilon$, this is $o(e^6)$ and does not
specify the coefficient of $e^8$.
The comparison here is structural. A coefficient identification
in physical momentum would also require the exact branch and conjugacy
identifications above and the evaluation of the intermediate diagonals.
For Corollary~\ref{st-cor:eighth-detuning} the additional diagonal combination
to be evaluated, with $\widetilde r_n=r_n+r_n^{(\infty)}$,
would be
\[
  2(\widetilde r_3+\widetilde r_{-3})
       -\widetilde r_{12}-\widetilde r_{-12}.
\]
The transport, dispersion, and sign terms cancel in this
combination. The physical combinations in \eqref{st-eq:main}
and \eqref{st-eq:eighth} are computed directly.
\end{remark}

\section{Consequences for frequency relations and dispersion laws}
\label{st-sec:consequences}
\subsection{Finite frequency relations}
\begin{corollary}\label{st-cor:relations}
Fix a finite list of branches $(j_\ell,\tau_\ell)$ with nonsquare
$j_\ell\ge3$ and fixed real coefficients $c_\ell$ satisfying
$\sum_\ell c_\ell j_\ell=\sum_\ell c_\ell\tau_\ell\sqrt{j_\ell}=0$.
On a common small momentum interval for this list,
\[
 \sum_\ell c_\ell\lambda_{j_\ell,\tau_\ell}(I)
   =2\Bigl(\sum_\ell c_\ell\Bigr) I^3+O(I^4).
\]
The remainder constant may depend on the list and its coefficients.
\end{corollary}
\begin{proof}
The constraints cancel the $v$ and $d_s$ terms in
\eqref{st-eq:variance-expansion}. The sum is $(\sum_\ell c_\ell)m_s(e)
+O(e^8)$, and $e^6=8I^3+O(I^4)$. When the coefficient sum is nonzero,
the sign at sufficiently small positive $I$ is the sign of that sum.
\end{proof}

Thus a relation with $\sum_\ell c_\ell\ne0$ is broken at order six,
while one with $\sum_\ell c_\ell=0$ persists through that order.
For instance, fix $p\ge2$ and a nonsquare $q\ge3$, both integers.
The flat frequencies of $(p^2q,+)$, $(q,-)$ and $(q,+)$
satisfy the primitive integer relation with coefficients
$(1,-p(p-1)/2,-p(p+1)/2)$. Their combination has leading term
$-2(p^2-1)I^3$, with remainder $O_{p,q}(I^4)$. The case $p=2$ is
\eqref{st-eq:main}.

\subsection{Best fits by common dispersion coefficients}
For the two frequencies at a fixed index,
$\lambda_{j,\tau}=p_j+\tau\sqrt j\,s_j$, where $p_j$ is their
mean and $s_j$ their normalized gap. A common three-coefficient fit
is exact precisely when $s_j$ is constant and $p_j$ is affine on
the fitted indices, provided both signs are included at each index.
For $j\ge3$ nonsquare, a common additive constant absorbs the
sixth-order offset in \eqref{st-eq:variance-expansion}. Including
$j=2$ produces a branch-dependent cubic correction in momentum.
For pairs of nonsquare indices above $2$, the first discrepancy
in the normalized gaps appears at quartic order.
For nonsquare $q\ge2$ write $\mathscr S_q=\{(4q,+1),(q,-1),(q,+1)\}$.

\begin{corollary}
\label{st-cor:mean-law}
For each such $q$, on a sufficiently small momentum interval,
\[
 \begin{split}
 &\inf_{\mathfrak v,\mathfrak d\in\mathbb R}
 \max_{(j,\tau)\in\mathscr S_q}
 \bigl|\lambda_{j,\tau}(I)+j\mathfrak v-\tau\mathfrak d\sqrt j\bigr|\\
 &\hspace{3em}=\frac{|\Delta_q(I)|}{5}.
 \end{split}
\]
The coefficient of $I^3$ is $6/5$ for $q\ge3$ and
$(158+32\sqrt2)/105$ for $q=2$.
\end{corollary}
\begin{proof}
Use the order $(4q,+1),(q,-1),(q,+1)$. The model vectors
span the two-dimensional plane annihilated by $c=(1,-1,-3)$.
For $x=(\lambda_{4q,+},\lambda_{q,-},\lambda_{q,+})$ and every vector
$y$ in that plane,
$|\Delta_q|=|c\cdot(x-y)|\le5\|x-y\|_\infty$.
Conversely, set $r=(\Delta_q/5)(1,-1,-1)$. Then
$c\cdot r=\Delta_q$, so $x-r$ is in the model plane and
$\|r\|_\infty=|\Delta_q|/5$. This gives equality. The final expansion
follows from Theorem~\ref{st-thm:detuning}.
\end{proof}

\begin{corollary}
\label{st-cor:three-law}
Take $\mathscr S=\mathscr S_2\cup\mathscr S_3$, joining the boundary
triple to the first triple protected through sixth order. Set
\[
 E_3(I)=\inf_{\mathfrak v,\mathfrak d,\mathfrak m\in\mathbb R}
 \max_{(j,\tau)\in\mathscr S}
 |\lambda_{j,\tau}(I)+j\mathfrak v-\tau\sqrt j\,\mathfrak d-\mathfrak m|.
\]
Then as $I\downarrow0$,
\[
 E_3(I)=\frac{8(9+16\sqrt3)}{21(11+4\sqrt6)}I^3+O(I^4).
\]
The infimum permits a new choice of all three coefficients at each $I$.
\end{corollary}
\begin{proof}
Order the branches as $(8,+),(2,-),(2,+),(12,+),(3,-),(3,+)$ and
denote by $L$ the space of model vectors with entries
$-j\mathfrak v+\tau\sqrt j\,\mathfrak d+\mathfrak m$.
Subtract the vector with coefficients $(v,d_s,m_s)$ from the physical
frequencies. By \eqref{st-eq:variance-expansion} and the two $j=2$
corrections in \eqref{st-eq:boundary}, the residual is
$(8/21)I^3 b+O(I^4)$ in $\ell^\infty$, where
$b=(0,1-2\sqrt2,1+2\sqrt2,0,0,0)$. Here the comparison expansion
at $j=2$ is justified in Section~\ref{st-sec:return}.
Distance to $L$ is homogeneous and $1$-Lipschitz, so it suffices to
compute $\operatorname{dist}_\infty(b,L)$.

Set $\rho=\sqrt{3/2}$, $u=(9+16\sqrt3)/(11+8\rho)$,
$t=2\sqrt2-u$ and $V=1+\rho t-u$. The model
$P_{j,\tau}=1-(j-2)V+(\tau\sqrt j/\sqrt2)t$ has residual
\[
 b-P=(5+(6\rho-2)t-6u,-u,u,u,2\rho t-u,-u).
\]
Squaring the endpoints gives $1414/1000<\sqrt2<1415/1000$,
$1732/1000<\sqrt3<1733/1000$ and $1224/1000<\rho<1225/1000$.
Substitution yields $7/4<u<9/5$, $1<t<11/10$ and $6/5<\rho<5/4$.
The first residual lies between $-3/5$ and $11/20$, and the fifth
lies between $3/5$ and $1$. Their absolute values are below $u$.
For the lower bound, $c=(0,8\rho-9,-8\rho-9,-2,0,20)$
annihilates the three model vectors and satisfies
$\|c\|_1=22+16\rho$ and $|c\cdot b|=18+32\sqrt3$.
Hence $\|b-P'\|_\infty\ge |c\cdot b|/\|c\|_1=u$ for every $P'\in L$.
The displayed fit attains this lower bound.
\end{proof}

\begin{corollary}
\label{st-cor:eighth-fit}
Fix distinct nonsquare $j,k\ge2$. On a sufficiently small common
momentum interval for their four physical branches, the best uniform error $E_{j,k}(I)$ of the three-coefficient model is
\[
 E_{j,k}(I)=\frac{\sqrt{jk}}{\sqrt j+\sqrt k}|s_j(I)-s_k(I)|.
\]
For $j,k\ge3$ the leading term is of order $I^4$ and follows from
\eqref{st-eq:gap-variation}. For $j=2$ and $k\ge3$,
$s_2-s_k=(16/21)I^3+O_k(I^4)$, so the leading error is of order
$I^3$. At $j=3$ and $k=12$, the exact error is $|Q(I)|/6$.
\end{corollary}
\begin{proof}
In the order $(j,+),(j,-),(k,+),(k,-)$, the three model
vectors span a hyperplane. Indeed, the two pairwise means separate
the constant and linear functions of the index, and the two
differences give a third direction. The normal to this hyperplane
is $c=(\sqrt k,-\sqrt k,-\sqrt j,\sqrt j)$.
If $x$ is the vector of the four physical frequencies, its
distance in $\ell^\infty$ to the hyperplane is
$|c\cdot x|/\|c\|_1$: H\"older's inequality gives the lower bound
and subtracting $(c\cdot x)\operatorname{sign}(c)/\|c\|_1$
attains it, with the sign taken componentwise.
Here $c\cdot x=2\sqrt{jk}(s_j-s_k)$ and
$\|c\|_1=2(\sqrt j+\sqrt k)$.
For index $2$, the two corrections in \eqref{st-eq:boundary} change
the normalized gap by $(2/21)e^6+O(e^8)$. The comparison gap
through order six has the same formula at all indices in question,
by the calculation in Section~\ref{st-sec:return}.
\end{proof}

\subsection{The mean of the two signs}
The pair mean $p_j=(\lambda_{j,+}+\lambda_{j,-})/2$ removes the
dispersion coefficient from any common three-coefficient fit.
The remaining two coefficients would require the $p_j$ to be
affine functions of the Fourier index.

\begin{corollary}
\label{st-cor:mean-curvature}
For any fixed nonsquare $3\le j<k<\ell$, set
\[
 D_{j,k,\ell}(I)=\frac{p_k(I)-p_j(I)}{k-j}
                      -\frac{p_\ell(I)-p_k(I)}{\ell-k}.
\]
Then
\[
 D_{j,k,\ell}(I)=16I^4
 \left(\frac{u_k-u_j}{k-j}-\frac{u_\ell-u_k}{\ell-k}\right)
 +O_{j,k,\ell}(I^5),
\]
where
\[
 u_n=c(n)+\frac{\mathbf 1_{\{3\}}(n)}{528},\qquad
 c(x)=\frac3{16}+\frac9{128(4x-1)}-\frac9{128(4x-9)} .
\]
The coefficient of $I^4$ is positive.
On a sufficiently small common momentum interval for $j,k,\ell$,
the best uniform error for the three means by an affine function
of the index is
\[
 \frac{(k-j)(\ell-k)}{2(\ell-j)}|D_{j,k,\ell}(I)|.
\]
This is also a lower bound for the best common three-coefficient
fit on the six original branches at $j,k,\ell$.
\end{corollary}
\begin{proof}
The even part of $g(w)$ in Section~\ref{st-sec:eighth} is $c(w^2)$.
Averaging the two $C_{3,\tau}$ gives $1/528$. Thus, with
$W_5=[e^5]W$, the fixed-index expansion is
\[
 p_n=-nv+e^6/4+e^8\bigl(\widehat W_5(1)/2+u_n\bigr)+O_n(e^{10}).
\]
The affine terms cancel in the difference of slopes, giving the
stated expansion for $D_{j,k,\ell}$.
For $x\ge5$, $c''(x)=(9/4)((4x-1)^{-3}-(4x-9)^{-3})<0$.
To compare a slope starting at $3$ with the next slope, use
\[
 c(5)-c(3)-1/528-2c'(5)=19645/2096688>0 .
\]
The function $c(x)-c(3)-1/528-(x-3)c'(x)$ has derivative
$-(x-3)c''(x)>0$ for $x\ge5$. Thus the slope from $3$ to $k$
is greater than $c'(k)$, which is greater than the slope from
$k$ to $\ell$. This proves strict positivity in both cases.

Let $L_{j,\ell}$ be the line through $(j,p_j)$ and $(\ell,p_\ell)$.
The middle value lies above the chord for sufficiently small positive
$I$, since its height relative to the chord is
\[
 H_k=p_k-L_{j,\ell}(k)=
       \frac{(k-j)(\ell-k)}{\ell-j}\,D_{j,k,\ell}.
\]
Subtracting any affine function from all three means leaves $H_k$
unchanged. If the residuals of that fit are bounded by $E$, the
residual at $k$ minus the interpolated endpoint residuals has
magnitude at most $2E$. The line $L_{j,\ell}+H_k/2$ attains
that bound. Averaging the errors at each pair of signs can only
decrease the maximum absolute error, giving the bound for six branches.
\end{proof}

\subsection{Scalar sublevel bounds and momentum scales}
The frequency combination $\Delta_q$ has a cubic zero in the
momentum. The following bounds measure its small-value sets.

\begin{corollary}\label{st-cor:sublevel}
For each fixed nonsquare $q\ge3$ there is $0<I_q^*<I_q$ such that
\begin{equation}\label{st-eq:detuning-derivatives}
 -24I^2\le\Delta_q'(I)\le-12I^2,\qquad
 -48\le\Delta_q'''(I)\le-24,
 \qquad 0\le I\le I_q^*.
\end{equation}
In particular $\Delta_q$ is strictly decreasing. For every fixed threshold
$a\in\mathbb R$ and window radius $\delta>0$,
\[
 \meas\{I\in[0,I_q^*]:|\Delta_q(I)-a|<\delta\}
 \le\min\{I_q^*,(\delta/2)^{1/3}\}.
\]
If $0<I_0$ and $2I_0\le I_q^*$, then also
\[
 \meas\{I\in[I_0,2I_0]:|\Delta_q(I)-a|<\delta\}
 \le\min\{I_0,\delta/(6I_0^2)\}.
\]
At $a=0$ and as $\delta\downarrow0$, the sublevel measure on
$[0,I_q^*]$ equals
$(\delta/6)^{1/3}+O_q(\delta^{2/3})$. For a moving threshold
$a:[I_0,2I_0]\to\mathbb R$ with Lipschitz constant at most $6I_0^2$,
the sublevel measure on $[I_0,2I_0]$ with $a$ replaced by $a(I)$
is at most $\min\{I_0,\delta/(3I_0^2)\}$.
For $q=2$, at zero threshold and on a sufficiently small fixed
momentum interval, the sublevel length is
$(\delta/c_2)^{1/3}+O(\delta^{2/3})$, with $c_2$ from
\eqref{st-eq:main}.
\end{corollary}
\begin{proof}
By the analytic assertion of Theorem~\ref{st-thm:detuning},
\[
 \Delta_q(I)=-6I^3+I^4r_q(I)
\]
with analytic $r_q$ near zero. On a sufficiently small compact interval
the functions $r_q$ and its first three derivatives are bounded.
Since
\[
 \Delta_q'(I)=-18I^2+I^3(4r_q+Ir_q'),
\]
we can choose $I_q^*$ so that
$I|4r_q+Ir_q'|\le6$ there. Differentiating the analytic remainder
three times gives
\[
 \Delta_q'''(I)=-36+24Ir_q+36I^2r_q'
                         +12I^3r_q''+I^4r_q'''.
\]
Reduce the same radius so the last four terms have absolute value
at most $12$. This proves \eqref{st-eq:detuning-derivatives}.
All choices are made after $q$ is fixed.

For $0\le u<v\le I_q^*$, integration gives
\[
 \Delta_q(u)-\Delta_q(v)\ge12\int_u^v I^2\,dI
 =4(v^3-u^3)\ge4(v-u)^3.
\]
The inverse image of $(a-\delta,a+\delta)$ is an interval, possibly
empty, because $\Delta_q$ is continuous and strictly decreasing.
For any two of its points the left side is less than $2\delta$;
its diameter, and hence its Lebesgue measure, is at most
$(\delta/2)^{1/3}$. The ambient interval gives the other global bound.
On $[I_0,2I_0]$, the estimate $\Delta_q(u)-\Delta_q(v)\ge12I_0^2(v-u)$
gives the stated bound on that interval by the same argument.

For $a=0$, the upper endpoint $I_\delta$ of the set solves
$-\Delta_q(I_\delta)=\delta$ for all sufficiently small $\delta>0$.
The sign bound gives $I_\delta=O_q(\delta^{1/3})$. Substituting this
into $\delta=6I_\delta^3(1+O_q(I_\delta))$ gives the stated expansion
and hence the sharp exponent $1/3$.
At $q=2$, analyticity and $c_2>0$ in \eqref{st-eq:main} give
strict monotonicity near zero and the same endpoint argument with
$6$ replaced by $c_2$.
For a moving $a$, the difference $(\Delta_q-a)(u)-(\Delta_q-a)(v)$
for $u<v$ is bounded below by $12I_0^2(v-u)-|a(v)-a(u)|\ge6I_0^2(v-u)$.
The same diameter argument gives the moving-threshold bound.
\end{proof}

\begin{corollary}\label{st-cor:scales}
There are constants $c,C,c_*>0$ independent of nonsquare $q\ge3$.
On $|I|\le c/q$ the analytic continuation satisfies
$|\Delta_q(I)+6I^3-b_q I^4|\le Cq^{3/2}|I|^5$, where the Taylor
coefficient $b_q=[I^4]\Delta_q$ satisfies $|b_q|\le C$.
On $0\le I\le c_*/q$ the cubic coefficient controls the sign and
the derivative estimates in \eqref{st-eq:detuning-derivatives} hold.
\end{corollary}
The proof is given in Appendix~\ref{st-sec:polynomial}.
The same spectral circles show that all the normalized gaps with
nonsquare indices in $[2,N]$, where $N$ is an integer Fourier
cutoff, are ordered on a common positive momentum interval of length
proportional to $N^{-3}$
(Corollary~\ref{st-cor:uniform-gap}).

\subsection{Coefficient checks in the original conformal coordinate}\label{st-sec:coefficient-check}
The coefficients can also be computed before the transport change of
coordinate. In \eqref{st-eq:A}, write $a=\sum_r e^ra_r$ and
$b=\sum_r e^rb_r$. At order $r>0$ the matrix from input $n$ to output
$m$ is
\[
 \begin{pmatrix}
 m\widehat b_r(m-n)&0\\
 \widehat a_r(m-n)&n\widehat b_r(m-n)
 \end{pmatrix}.
\]
Generate $a_r,b_r$ from \eqref{st-eq:recursion} and the definitions
following \eqref{st-eq:velocities}. At Taylor order $N$, the speed
coefficient at order $N$ is fixed by projecting the Babenko equation
at order $N+1$ onto the first harmonic.
The recursion \eqref{st-eq:right-recursion}, now applied directly to
$A_e$, retains both components of every index in $[j-r,j+r]$ at order
$r$. An order-$s$ coefficient shifts the index by at most $s$, and the
scalar eigenvalue feedback preserves the index. These intervals
therefore include all contributions to the desired jet, including zero
and negative intermediate indices. The reduced solve at $j$ retains
the other eigenvector of its $2\times2$ block.

The ancillary programs \texttt{check-conformal-jets.py} and
\texttt{check-sharp-protection.py} implement these recurrences
over the rationals and $\mathbb Q(\sqrt j)$. The \texttt{README.md}
in the \texttt{checks} directory gives the pinned dependencies,
reproduction commands, and archived outputs. From that directory,
after installing the listed dependencies, run
\texttt{python check-conformal-jets.py --order 8 --output fresh.json}
to write a fresh result file. In the table, $[e^r]\Delta_q$ means the
coefficient after composing $\Delta_q$ with $I(e)=\mathcal I(U_e)$.
\begin{center}
\begin{tabular}{@{}crrrr@{}}
\toprule
$q$ & $[e^2]\Delta_q$ & $[e^4]\Delta_q$ & $[e^6]\Delta_q$ & $[I^3]\Delta_q$ \\
\midrule
$2$ & $0$ & $0$ & $-(79+16\sqrt2)/84$ & $-(158+32\sqrt2)/21$ \\
$3,5,6,7,12$ & $0$ & $0$ & $-3/4$ & $-6$ \\
\bottomrule
\end{tabular}
\end{center}
For the four-branch relation in Corollary~\ref{st-cor:eighth-detuning}, the
program evaluates the two conformal recurrences with upper-right
symbols $|n|$ and $n$ separately. It verifies vanishing through
amplitude degree seven and gives the following entries:
\[
\begin{array}{c|ccc}
 & \text{physical} & \text{comparison} & \text{difference}\\ \hline
[e^8]Q /\sqrt3 & 1102/20163 & 1513/53768 & 7/264
\end{array}
\]
The sixth-order coefficient also vanishes in $Q(I(e))$, so the
momentum coefficient follows from $e^8=16I^4+O(I^5)$.
The second script also compares upper-right
symbols $|D|$ and $D$ in the original conformal coordinates and
checks \eqref{st-eq:sharp-protection} at $j=2,3,5,6$, for both signs.
The identity $C_F^\sharp D C_F=D$ in \eqref{st-eq:congruence-D}
identifies the $D$ implementation with the comparison operator.
The scripts check finite Taylor jets. The preceding proofs establish
the analytic continuation and remainder bounds for the selected
branches of the full operators.

The sublevel estimates here concern scalar frequency combinations on
the Stokes branch. The almost-periodic existence question
in~\cite[Section~1.4, P1,3]{BertiSurvey} also involves coupled normal
blocks and nonlinear iteration. For the frequencies considered here,
the explicit leading coefficients determine the momentum scales on
which the selected flat relations split.

\appendix
\section{Analytic realization of the conformal operators}
\label{st-sec:conformal}

\subsection{The spatial analytic estimate}
For $\rho>0$, let $\mathcal A_\rho$ be the Fourier--Wiener space
with norm
\[
 \|f\|_\rho=\sum_{j\in\mathbb Z}e^{\rho|j|}|\widehat f(j)|.
\]
It is a Banach algebra. We first estimate the Fourier kernel for a
small circle map on a complex analytic ball. The energy identity
then identifies that kernel for the physical boundary map on the real
slice and gives the maximal-domain statement.

\begin{theorem}
\label{st-thm:dno}
Let $0<\rho_0\le1$, put $\rho=\rho_0/8$, and suppose that a real surface
$\eta\in\mathcal A_{\rho_0}$ satisfies
$\|\eta\|_{\rho_0}\le\rho_0/1024$. Then the infinite-depth physical
Dirichlet--Neumann operator has the form
\[
 G(\eta)=|D|+R_\eta,\qquad R_\eta1=0,\qquad
 \langle R_\eta f\rangle=0 .
\]
Its matrix coefficients, defined by
$(R_\eta)_{mj}=\widehat{R_\eta e^{ijx}}(m)$, satisfy
\begin{equation}\label{st-eq:dno-kernel}
 |(R_\eta)_{mj}|
 \le C\rho^{-3}\|\eta\|_{\rho_0}
                         e^{-\rho(|m|+|j|)/4}.
\end{equation}
For nonnegative integers $p,r$, define
\[
 \|R\|_{\mathfrak K_{p,r}}
 =\sum_{m,j}\langle m\rangle^p\langle j\rangle^r
        e^{\rho(|m|+|j|)/8}|R_{mj}|.
\]
Then
\[
 \|R_\eta\|_{\mathfrak K_{p,r}}
       \le C_{p,r}\rho^{-p-r-5}\|\eta\|_{\rho_0}.
\]
The map $\eta\mapsto R_\eta$ has a holomorphic extension on the
complex open ball $\|\eta\|_{\rho_0}<\rho_0/512$, with values in
$\mathcal B(L^2,L^2)$.
For real surfaces the distributional maximal domain of $G(\eta)$
on $L^2$ is $H^1$, and its restriction to $L^2_0$ is self-adjoint,
nonnegative, with domain $H^1_0$.
\end{theorem}

The construction allows nonzero surface means: vertical translation
does not change $G(\eta)$ in infinite depth. The main text restricts to
zero physical mean. In particular the holomorphic extension is
defined on an ordinary open function ball before any constraint is
imposed.

\subsection{Circle maps and the two Hardy sectors}
We first treat a general small circle map $F(y)=y+h(y)$, not
necessarily a conformal boundary map. This applies both to the
transport map in Section~\ref{st-sec:stokes} and to the conformal map
constructed below. Set $C_Fu=u\circ F$ and assume
\[
 h\in\mathcal A_{\rho_0/2},\qquad
 \|h\|_{\rho_0/2}\le\rho_0/1024.
\]
For $r>0$, the Fourier series and Banach algebra estimate give
$\|u\circ(\mathrm{Id}+k)\|_r\le\|u\|_{r+\|k\|_r}$.
A derivative costs at most $1/\epsilon$ when shrinking a strip by
$\epsilon>0$, by $|n|e^{-\epsilon|n|}\le1/\epsilon$.
For the inverse $F^{-1}=\mathrm{Id}+k$, the map
$k\mapsto-h\circ(\mathrm{Id}+k)$ takes the ball
$\|k\|_{\rho_0/4}\le\rho_0/16$ into itself. Along segments in this
ball there remains a strip of width $3\rho_0/16$ for differentiating $h$,
so the contraction constant is at most
$16\|h\|_{\rho_0/2}/(3\rho_0)<1/4$.
This also applies on the larger complex ball
$\|h\|_{\rho_0/2}<\rho_0/512$. For real $h$, $F'>1/2$.
On the real slice the $L^2$ adjoint is
\[
 C_F^\sharp=C_F^{-1}M_{1/F'}.
\]
On analytic test functions with a strip reserve for composition,
the same ordered expression defines $\sharp$ for complex $h$.
For the spectral arguments we use the bounded continuation of the
remainder constructed below.

In this calculation $x=F(y)$ denotes the image coordinate.
For real maps the chain rule and change of variables give the exact identities
\begin{equation}\label{st-eq:congruence-D}
 D_yC_F=M_{F'}C_FD_x,\qquad C_F^\sharp D_yC_F=D_x.
\end{equation}
We identify the two Fourier spaces when no coordinate subscript
is needed. On Fourier polynomials with $h$ real, set
$R_F=C_F^\sharp|D|C_F-|D|$. The matrix of $C_F$ has entries from input
$j$ to output $n$ given by
\[
 A_{nj}=\frac1{2\pi}\int_0^{2\pi}e^{ijF(y)-iny}\,dy.
\]
For real maps, change of variables gives
$(C_F^\sharp)_{mn}=A_{-n,-m}=\overline{A_{nm}}$.
The adjoint entry is continued to complex $h$ by $A_{-n,-m}$.
For real $h$ and positive input index $j$,
\begin{equation}\label{st-eq:hardy}
 (R_F)_{mj}=2\sum_{n<0}|n|A_{-n,-m}A_{nj},
\end{equation}
whereas for negative input index $j$,
\[
 (R_F)_{mj}=2\sum_{n>0}n A_{-n,-m}A_{nj}.
\]
These follow by subtracting, respectively, the two identities
$C_F^\sharp D C_F=D$ and $-C_F^\sharp D C_F=-D$.
Constants are preserved by $C_F$, and $|D|$ kills constants.
The real adjoint identity then gives both $R_F1=0$ and
$\langle R_Ff\rangle=0$.
For complex $h$ the same two series, with zero constant row and
column, define the continued matrix once the convergence below
is established.

\begin{lemma}\label{st-lem:kernel}
For $\|h\|_{\rho_0/2}\le\rho_0/1024$, the estimate
\eqref{st-eq:dno-kernel} holds for $R_F$ with
$\|\eta\|_{\rho_0}$ replaced by $\|h\|_{\rho_0/2}$.
Its Fr\'echet derivative in direction $u\in\mathcal A_{\rho_0/2}$
obeys
\[
 |D_h(R_F)_{mj}[u]|
 \le C\rho^{-3}\|u\|_{\rho_0/2}
                         e^{-\rho(|m|+|j|)/4}.
\]
The map $h\mapsto R_F$ is holomorphic into $\mathcal B(L^2,L^2)$
on the open ball $\|h\|_{\rho_0/2}<\rho_0/512$.
\end{lemma}
\begin{proof}
On the contours used below, $|h|\le H_0\le\rho/128$.
Moving the integration contour by
$\operatorname{sign}(j-n)\rho$ gives
\[
 |A_{nj}|\le e^{-\rho|n-j|+H_0|j|}.
\]
For $j>0$, write $J=j$ and $M=|m|$. Equation~\eqref{st-eq:hardy}
then gives
\[
 |(R_F)_{mj}|
 \le2e^{-\rho J+H_0(J+M)}
       \sum_{N\ge1}N e^{-\rho(N+|-N-m|)}.
\]
If $m\ge0$, the sum is at most $C\rho^{-2}e^{-\rho M}$.
If $m<0$, its exact value is
\[
 e^{-\rho M}\left\{\frac{M(M+1)}2+
        \frac{M}{e^{2\rho}-1}
        +\frac{e^{-2\rho}}{(1-e^{-2\rho})^2}\right\}
 \le C\rho^{-2}(1+\rho M)^2e^{-\rho M}.
\]
This includes the mixed-sign block, where just one factor has
opposite input and output signs. Absorbing the quadratic weight
in the exponential proves
\[
 |(R_F)_{mj}|\le C\rho^{-2}e^{-\rho(M+J)/2}.
\]
The formula for negative input indices gives the same estimate.
The zero row and column vanish exactly.

Differentiation of a composition entry gives
\[
 D_hA_{nj}[u]=\frac{ij}{2\pi}
                    \int_0^{2\pi}u\,e^{ijF-iny}\,dy .
\]
Differentiating both factors in \eqref{st-eq:hardy} therefore
adds at most $(M+J)\|u\|_{\rho_0/2}$ to the same convergent
majorant. Absorbing $M+J$ into the exponential gives the claimed
derivative bound. Integrating along $sh$, $0\le s\le1$,
and using $R_{\mathrm{Id}}=0$ gives the small value estimate.
The weighted double sum costs $C_{p,r}\rho^{-p-r-2}$,
which gives the stated kernel norm.

The integral entries are holomorphic in the ordinary complex
function ball. All sums and their first derivatives converge
locally uniformly there, with the same estimates, allowing
$H_0\le\rho/64$ on the larger open ball.
For fixed nonnegative integers $p,r$, the matrix-unit operator
from the input mode $j$ to the output mode $m$ has norm
$\langle m\rangle^p\langle j\rangle^r$ from $H^{-r}$ to $H^p$.
The exponential bound therefore gives locally uniform convergence
of finite matrix cutoffs in $\mathcal B(H^{-r},H^p)$, in particular
in $\mathcal B(L^2,L^2)$. Their holomorphic limits give the
operator-valued continuation.
\end{proof}

\begin{remark}\label{st-rem:H-not-small}
Let $P_-$ project onto the strictly negative Fourier indices.
Then
\[
 \mathcal H_F=C_F^\sharp(|D|-D)C_F
                    =(|D|-D)+R_F .
\]
Its full flat part $|D|-D=2|D|P_-$ has order one. On positive
inputs its entries equal those of $R_F$.
This is the distinction used in the full-space spectral homotopy.
\end{remark}

\subsection{Identification with the physical boundary problem}
For a small surface $\eta$, solve the conformal equation
\[
 h=\mathcal H[\eta\circ(\mathrm{Id}+h)],\qquad
 X_\eta=\mathrm{Id}+h,\qquad Y=\eta\circ X_\eta .
\]
Here $X_\eta$ is the map denoted by $X$ on the Stokes branch
in the main text. On $\|h\|_{\rho_0/2}\le\rho_0/8$,
composition uses at most the width $5\rho_0/8$ of $\eta$.
The remaining width bounds its derivative. Since the Hilbert
transform has norm at most one on $\mathcal A_{\rho_0/2}$,
the contraction constant is bounded by
$8\|\eta\|_{\rho_0}/(3\rho_0)<1/4$, and the fixed point satisfies
\[
 \|h\|_{\rho_0/2}\le\|\eta\|_{\rho_0}.
\]
On the open ball $\|\eta\|_{\rho_0}<\rho_0/512$, iterating
the contraction from zero gives holomorphic maps of $\eta$ by the
composition estimate on the reserved strip. The iterates converge
locally uniformly in $\mathcal A_{\rho_0/2}$, giving a holomorphic
map $\eta\mapsto h$. The bound on $h$ keeps it in the complex
ball for the continued kernel above.

For real $\eta$, $h=\mathcal HY$, so $h+iY$ has only nonpositive
Fourier frequencies. Its holomorphic extension to the lower
half-plane defines
\[
 Z(\zeta)=\zeta+h(\zeta)+iY(\zeta),\qquad
 |Z'-1|\le4\|\eta\|_{\rho_0}/\rho_0<1/2.
\]
Integration along segments of this convex half-plane shows
$|Z(\zeta_1)-Z(\zeta_2)|\ge|\zeta_1-\zeta_2|/2$.
Thus $Z$ is univalent and its boundary is
$(X_\eta(y),\eta(X_\eta(y)))$.
The displacement $Z-\mathrm{Id}$ is bounded at infinite depth.
The image lies below the graph: a crossing would contradict
injectivity on the closed half-plane. It is relatively open and
closed in that graph domain. Indeed preimages of a convergent
image sequence with interior limit cannot escape to infinity,
by bounded displacement, or approach the boundary.
Connectedness gives the entire physical fluid domain.

The pullback of the harmonic problem is consequently the Poisson
extension in the lower half-plane. Invariance of Dirichlet energy
and polarization give
\begin{equation}\label{st-eq:energy}
 \langle\psi_1G(\eta)\psi_2\rangle_x
     =\langle C_{X_\eta}\psi_1\,|D|C_{X_\eta}\psi_2\rangle_y,
 \qquad G(\eta)=C_{X_\eta}^\sharp|D|C_{X_\eta}.
\end{equation}
The displayed pairing is bilinear. For a real surface its
sesquilinear version follows by taking $\psi_1=\overline f$ and
$\psi_2=g$ for arbitrary test functions $f,g$, since composition
with $X_\eta$ commutes with conjugation.
Lemma~\ref{st-lem:kernel} now proves the kernel assertions of
Theorem~\ref{st-thm:dno} for the physical operator $G(\eta)$.
The composition $\eta\mapsto h\mapsto R_{X_\eta}$ proves its
holomorphic extension.

\subsection{The potential quotient and maximal domains}
Composition fixes constants but need not preserve mean zero.
On zero-mean representatives of the potential quotient set
\[
 C_0=P_0C_F|_{H^s_0},\qquad
 C_0^{-1}=P_0C_F^{-1}|_{H^s_0}.
\]
These are inverse maps because a discarded constant remains
constant under composition. For a real map $F$, $C_0^\sharp$
denotes the $L^2_0$ adjoint. Since $|D|P_0=P_0|D|=|D|$,
\[
 P_0C_F^\sharp|D|C_FP_0=C_0^\sharp|D|C_0 .
\]
The finite-rank mean repairs are thus part of the operator identity.

For a real diffeomorphism $F$ as above, change of variables bounds composition
and its adjoint on $L^2$. The chain rule gives positive integer
Sobolev bounds, interpolation gives intermediate orders, and
duality gives negative orders.
The energy identity \eqref{st-eq:energy} extends by density to
$\psi_2\in L^2$ and $\psi_1\in H^1$, with the left side interpreted
as a distributional pairing. This uses $|D|:L^2\to H^{-1}$ and
boundedness of $C_{X_\eta}$ at both test and data regularities.
Thus $G=|D|+R_\eta$ holds as a map $L^2\to H^{-1}$.
As $R_\eta$ is bounded on $L^2$,
\[
 u\in L^2,\quad G(\eta)u\in L^2
               \quad\Longleftrightarrow\quad u\in H^1.
\]
For real $\eta$, $R_\eta$ is a bounded self-adjoint perturbation
of $|D|$; nonnegativity follows from Dirichlet energy.
Fourier polynomials form a graph core, and restriction to
$L^2_0$ gives the domain $H^1_0$.
For complexified surface profiles $\eta$, the analytic
continuation is realized on the same domain $H^1$ in $L^2$.

\section{Arithmetic isolation and localized remainder bounds}
\label{st-sec:polynomial}

For a nonsquare integer $j\ge2$, put
$\delta_j=\operatorname{dist}(2\sqrt j,\mathbb Z)$.
This arithmetic distance gives a lower bound, up to a constant,
for the gap to flat frequencies of the opposite square-root sign.
We use it to choose the isolating contour. To estimate the correction inside that contour,
we use the Fourier localization of the factors adjoining the
exact mean inverse.

\begin{proposition}[Analytic discs and uniform eighth-order expansion]
\label{st-prop:polynomial}
There are constants $c_0,c,C,c_*>0$ independent of indices and signs.
For nonsquare $j\ge2$ the physical and comparison operators have
simple analytic branches on $|e|<2a_j$. The radius $r_j$ of the isolating
spectral contour and the amplitude scale $a_j$ are
\[
 r_j=\delta_j/32,\qquad a_j=c_0(\delta_j/\sqrt j)^{1/2}.
\]
For nonsquare $j\ge3$, $w=\tau\sqrt j$ and $|e|\le a_j/2$,
\begin{equation}\label{st-eq:polynomial-individual}
\begin{aligned}
 \left|\lambda_{j,\tau}+jv-w\sqrt M
           -\frac{1-w}{4}e^6-B_{j,\tau}e^8\right|
 &\le Cj^{3/2}|e|^{10},\\
 B_{j,\tau}&=\frac{1-w}{2}\widehat W_5(1)+g(w)
                  +\mathbf 1_{\{3\}}(j)C_{j,\tau}.
\end{aligned}
\end{equation}
Here $W_5=[e^5]W$, and $g$ is given in Section~\ref{st-sec:eighth}.
For nonsquare $q\ge3$ there is a real $b_q$, $|b_q|\le C$, such
that on $|I|\le c/q$,
\begin{equation}\label{st-eq:polynomial-momentum}
 |\Delta_q(I)+6I^3-b_q I^4|\le Cq^{3/2}|I|^5.
\end{equation}
The sign and derivative bounds of Theorem~\ref{st-thm:detuning}
and Corollary~\ref{st-cor:sublevel} hold on $0\le I\le c_*/q$.
\end{proposition}

\begin{proof}
\emph{The flat circle.}
Fix a complex amplitude neighborhood using the strips of
Lemma~\ref{st-lem:transported-type-a}, before selecting $j$.
Put $\beta_0=j-\tau\sqrt j$ and let $A_{\alpha,0}=A_\alpha(0)$
for the family \eqref{st-eq:homotopy}. For $n>0$, put
$h=n-j$ and $S=\sqrt n+\sqrt j$. The same-sign difference
for $n\ne j$ has modulus $|h|\,|1-\tau/S|\ge |h|/2$.
The opposite-sign difference is
\[
 |h+\tau S|=\frac{S}{S-\tau}\,
                 |h-(1-2\tau\sqrt j)|\ge\delta_j/2 .
\]
The nearest integer $m$ to $2\sqrt j$ satisfies
$|m^2-4j|\ge1$, so $\delta_j\ge1/(5\sqrt j)$.
Negative indices and the zero block have spectral real part at most
zero, while $\beta_0>0$. On the circle $|\zeta-\beta_0|=r_j$,
$|\zeta|\le2j$ and $\operatorname{Re}\zeta\ge j/4$.

For $j/2\le n\le8j$, each factor in
$(n-\zeta)^2-n=(n-\sqrt n-\zeta)(n+\sqrt n-\zeta)$ is bounded below by $r_j$,
and one has modulus at least $\sqrt n$. The inverse numerator is
$O(j)$, hence these blocks are bounded by $C\sqrt j/\delta_j$.
For $1\le n<j/2$ and large $j$, both factors are bounded below by
$j/4$. The finitely many smaller $j$ are covered by the spectral
separation above. For $n>8j$ the inverse is $O(1/n)$.
For $n=-m<0$, uniformly in $0\le\alpha\le1$,
\[
 |(m+\zeta)^2-m(2\alpha-1)|\ge(m+j/4)^2-m\ge(m+j)^2/16 .
\]
The last inequality reduces to $m(15m+6j-16)\ge0$.
This also covers the nondiagonalizable block at $\alpha=1/2$.
The triangular zero block has an inverse bounded on the same contour.
Thus the flat resolvent maps $L^2$ into $H^1$ and its $L^2$ norm
is at most $C\sqrt j/\delta_j$.

Set $S_e=\operatorname{diag}(1,v(e))$ and normalize the spectral
operator by $\widehat A_\alpha=v^{-1}S_e^{-1}A_\alpha S_e$.
It has diagonal $D$, upper-right
$D+\alpha(|D|-D+R_F)$, and lower-left $\kappa/v^2$.
The difference from $A_{\alpha,0}$ is bounded on $L^2$ by
$C|e|^2$ by \eqref{st-eq:kap3} and
Lemma~\ref{st-lem:quadratic-transport}.
Choose $c_0$ small enough that on $|e|<2a_j$ the product with
the flat resolvent has norm at most $1/2$. This gives a rank-one
Riesz branch $\widehat\beta_\alpha$ inside the contour, uniformly
in $\alpha$. The closed domain is $H^1$. On a possibly smaller
high-Sobolev disc this is the same branch as that used to extract
Taylor coefficients, by uniqueness of the Riesz branch.
The signed frequency is $\lambda_\alpha=-v\widehat\beta_\alpha$.

\emph{Two localized factors of the mean inverse.}
For the expansion we first use the differential endpoint $\alpha=0$.
Write $\widehat A_{\rightarrow}=\widehat A_\alpha|_{\alpha=0}$.
The similarity $\mathcal S_j=\operatorname{diag}(j^{-1/2},1)$
preserves $H^1$ for each $j$, and in the balanced coordinates
$\mathcal S_j\widehat A_{\rightarrow}\mathcal S_j^{-1}=B+V$.
Write $\mu=M/v^2$. The exact mean operator $B$ in the balanced
coordinates has Fourier blocks
\[
 B(n)=\begin{pmatrix}n&n/\sqrt j\\ \sqrt j\,\mu&n\end{pmatrix}.
\]
At $n=j$, let $r=(-\tau/\sqrt\mu,1)^{\mathsf T}$ and
$l=(-\tau\sqrt\mu/2,1/2)$ be the right and left mean eigenvectors,
with $lr=1$, extended by the mode $j$ and its coefficient functional.
Let $P=rl$, $P_c=1-P$, and $\bar\beta=j-w\sqrt\mu$, using the analytic
square root near $1$. The perturbation is $V=V[W]$, where
$V[f]$ has only the lower-left entry $\sqrt j\,f/v^2$.
Thus $lVr=0$.

For $|e|\le a_j$ and $|\zeta-\beta_0|\le r_j$, write
$R_c(\zeta):P_c L^2\to P_c H^1$ for the inverse of $B-\zeta$
on the complement, extended by zero on $P L^2$.
Let $\Pi_j$ project onto the Fourier indices $|n-j|\le\sqrt j/4$.
The operator norms in this calculation are the $L^2$ norms in
the balanced coordinates. In these norms,
\[
 \|R_c(\zeta)\|\le C/\delta_j,\qquad
 \|R_c(\zeta)\Pi_j\|\le C .
\]
For positive $j/2\le n\le8j$ the two spectral projections of $B(n)$
are uniformly bounded in the balanced coordinates, as follows
from $B(n)-n$ divided by $\sqrt{n\mu}$.
In this range each root moves by at most $C\sqrt j |e|^2$, less than
$r_j$ after reducing $c_0$. In particular $|\bar\beta-\beta_0|\le r_j$.
At $n=j$ the selected branch is removed from $R_c(\zeta)$, and the
other root is at distance of order $\sqrt j$.
For positive $n$ outside that range, use
$|\sqrt{n\mu}|\le C\sqrt n$ in the two root distances.
For $n>8j$ both distances are bounded below by a multiple of $n$.
For $n<j/2$ and large $j$ they are bounded below by a multiple of $j$.
The finitely many remaining indices are controlled by the flat
separation and continuity.
For $n=-m<0$, the determinant $(m+\zeta)^2+m\mu$ has the same
lower bound up to a constant, since the change from $\mu=1$ is
$O(m|e|^2)$ and $m/(m+j)^2\le1/8$. The zero block is triangular.
These inverse numerators are $O(|n|+j)$ and decay of the inverse
at infinity gives $R_c(\zeta):L^2\to H^1$ at each fixed $j$.
For $0<|n-j|\le\sqrt j/4$, the same-sign root is at least
$|n-j|/2$ from $\bar\beta$, and the opposite-sign root is at least
$c\sqrt j$ away. This follows from $\sqrt\mu=1+O(e^2)$ and
$n-j=(\sqrt n-\sqrt j)(\sqrt n+\sqrt j)$. Since
$|\zeta-\bar\beta|\le 2r_j$, these estimates give the uniform bound
for $R_c(\zeta)\Pi_j$.

On a fixed spatial Wiener strip of width $\chi>0$, let $\|f\|_\chi$
be the norm from Appendix~\ref{st-sec:conformal}. A factor $V[f]$
attached to the selected mode has a tail outside $\Pi_j$ of size at most
$C\sqrt j\,e^{-\chi\sqrt j/4}\|f\|_\chi$ on either side.
Using $\delta_j^{-1}\le 5\sqrt j$ gives
\begin{equation}\label{st-eq:strip-wings}
 \|R_c(\zeta)V[f]r\|_2+\|lV[f]R_c(\zeta)\|_{2\to\mathbb C}
 \le C\sqrt j\|f\|_\chi .
\end{equation}
The projector is used on the column $V[f]r$ and the row $lV[f]$.
The exact mean inverse preserves Fourier indices, and the
exponential tails pay its global factor $\delta_j^{-1}$.
The full perturbed inverse below uses its global bound.

\emph{The remainder after the quadratic pairing.}
Set $\zeta=\widehat\beta_{\rightarrow}(e)$, the selected eigenvalue
of $\widehat A_{\rightarrow}$.
Since $W=O(e^3)$ on the fixed strip, $\|R_c(\zeta)V\|\le
C\delta_j^{-1}\sqrt j |e|^3$. For $|e|\le a_j$,
\[
 \sqrt j\,|e|^2/\delta_j\le c_0^2.
\]
Thus the Neumann norm is at most $Cc_0^2|e|<1/2$ after reducing $c_0$.
Write the eigenvector as $s r+u_c$ with $u_c\in P_cH^1$.
The complementary equation is
$(1+R_c(\zeta)V)u_c=-sR_c(\zeta)Vr$, so invertibility on $L^2$
forces $s\ne0$. After normalizing $s=1$, applying $l$ to the
eigenvalue equation gives
\[
 |\zeta-\bar\beta|
 =|lV(1+R_c(\zeta)V)^{-1}R_c(\zeta)Vr|
 \le Cj|e|^6 .
\]
The complementary solution lies in $H^1$ by substitution, since
$R_c(\zeta)$ takes $L^2$ into $H^1$. In particular the Neumann series
is used on $L^2$ with no implied $H^1$ operator bound uniform in $j$.

Now take $R=R_c(\bar\beta)$ and $d=\zeta-\bar\beta$.
On $P_c L^2$ put $E=P_c V P_c-dP_c$ and let $\mathcal G$ be the inverse
of $B-\bar\beta+E$ on this subspace, extended by zero.
Since $\|P\|+\|P_c\|$ is uniformly bounded,
$\|E\|\le C\sqrt j|e|^3$, and $\|\mathcal G\|\le C/\delta_j$
by the same Neumann estimate with a reduced $c_0$.
The resolvent identity $\mathcal G=R-RER+RE\mathcal G ER$ gives
\[
 d=-lVRVr+lVRERVr-lVRE\mathcal G ERVr .
\]
The two factors in \eqref{st-eq:strip-wings} bound the last term by
\[
 Cj^2\delta_j^{-1}|e|^{12}
 = Cj^{3/2}|e|^{10}\frac{\sqrt j\,|e|^2}{\delta_j}
 \le Cc_0^2 j^{3/2}|e|^{10}.
\]
Since $R=RP_c=P_cR$, the middle term splits into
$lVRVRVr-d\,lVR^2Vr$. The two localized factors and the bound
$|d|\le Cj|e|^6$ give $|d\,lVR^2Vr|\le Cj^2|e|^{12}$,
which obeys the same bound since $\delta_j\le1/2$.
In the remaining cubic term replace
$W$ by $e^3 W_3$ with $W_3=2\cos z$ in the spatial coordinate of
Section~\ref{st-sec:stokes}. Three shifts by $+1$ or $-1$ cannot
return to $j$, so the resulting term is zero even with the exact
mean in $R$. Errors in this replacement are bounded by
$Cj^{3/2}|e|^{10}$, using \eqref{st-eq:strip-wings} on both
ends and $\|W-e^3W_3\|_\chi\le C|e|^4$.
Consequently the quadratic pairing in \eqref{st-eq:eighth-schur},
after restoring the factor $v$ and the sign of $\lambda$,
has error $O(j^{3/2}|e|^{10})$ uniformly.

For that quadratic pairing, Taylor-expand $W$ through degree six
on the fixed strip. Terms with the $O(e^7)$ remainder have bound
$Cj|e|^{10}$, again by \eqref{st-eq:strip-wings}. The expansion
has spatial shifts $|h|\le4$ by Lemma~\ref{st-lem:highest}.
For each of these finitely many nonzero offsets the scalar kernel
in \eqref{st-eq:eighth-schur} is analytic and bounded by $Cj$
on a fixed small disc, independently of $j$.
For large $j$ this follows from the $2vw\sqrt M$ term in the
denominator, and for the remaining nonsquares it follows from
$h+2w-1\ne0$ at zero. Its Taylor coefficients through order
eight in the pairing are those computed in
Section~\ref{st-sec:eighth}. Half-period parity excludes odd
degrees. Cauchy's estimate on the fixed small disc bounds the
remaining terms of this finite sum by $Cj|e|^{10}$.

For nonsquare $j\ge5$, transfer this expansion to the physical
branch using Corollary~\ref{st-cor:uniform-comparison} below.
On $|e|\le a_j/2$, the comparison difference is at most
\[
 C\delta_j a_j^{-10}2^{-(2j-8)}|e|^{10}\le C|e|^{10},
\]
since $a_j^{-1}\le C\sqrt j$ and $\delta_j\le1/2$.
For $j=3$ the sharp comparison supplies $C_{3,\tau}e^8$
with error $O(e^{10})$. This proves \eqref{st-eq:polynomial-individual}.

\emph{Momentum and derivatives.}
In the detuning combination, the transport and dispersion
means cancel identically. The sixth-order term is $-3e^6/4$.
The combination of the $B_{j,\tau}$ is uniformly bounded:
the part proportional to $w$ cancels, and
$g(w)=-w/16+3/16+O(1/|w|)$ by division of its rational
expression. The correction at $j=3$ concerns a fixed index.
By \eqref{st-eq:action}, $e^2=S(I)=2I+O(I^2)$ on a fixed
complex disc. Since $a_j^2\ge c_0^2/(5j)$, a single choice
of $c$ ensures $|S(I)|\le a_j^2/4$ for $j=q,4q$ on
$|I|\le c/q$. All branches are even in $e$, so their compositions
with $S(I)$ through the variable $e^2$ are holomorphic on this
momentum disc. This gives \eqref{st-eq:polynomial-momentum}
with $b_q$ uniformly bounded.
Write $\Delta_q=-6I^3+b_q I^4+I^5\mathcal R_q(I)$.
On an inner disc Cauchy's estimate gives
$|\mathcal R_q^{(m)}|\le C_m q^{3/2+m}$ for $0\le m\le3$.
Thus the deviations in the first and third derivatives are
bounded respectively by $C(I^3+q^{3/2}I^4)$
and $C(I+q^{3/2}I^2)$ for $I\ge0$ on that disc.
The same bounds at derivative order zero control the sign.
Choosing $c_*$ small uniformly gives the claimed interval.
One may take $I_q=2c_*/q$ and $I_q^*=c_*/q$.
\end{proof}

\begin{corollary}
\label{st-cor:uniform-comparison}
Let $a_j$ be as above, $h_j=\lambda_{j,\tau}-\lambda_{\rightarrow,j,\tau}$
and $m_j=2j+2$. On $|e|\le a_j$,
\[
 |h_j(e)|\le C\delta_j(|e|/a_j)^{m_j},\qquad
 |h_j(e)-C_{j,\tau}e^{m_j}|
           \le 2C\delta_j(|e|/a_j)^{m_j+2}.
\]
\end{corollary}
\begin{proof}
The two rescaled eigenvalues lie in the same disc bounded by the
contour of radius $r_j=\delta_j/32$. Undo the common spectral rescaling by $v$.
Apply the maximum modulus principle to $h_j(e)/e^{m_j}$,
using Theorem~\ref{st-thm:sharp-protection} at zero.
Cauchy's formula also bounds $|C_{j,\tau}|a_j^{m_j}$ by
$C\delta_j$. By evenness the subtracted difference has a zero
of order at least $m_j+2$, which gives the second estimate
by the same principle.
\end{proof}

\begin{corollary}
\label{st-cor:uniform-gap}
There are constants $c,C,c_*>0$ independent of the indices.
For nonsquare $3\le j<k$ and $|I|\le c/k$,
\[
 \bigl|s_j-s_k-16(\mathfrak f(j)-\mathfrak f(k))I^4\bigr|
       \le C k |I|^5.
\]
For real $0<I\le c_*(k-j)/(j k^2)$ one has $s_j(I)>s_k(I)$.
In particular the nonsquare indices in $[3,N]$ are strictly
ordered on the common interval $0<I\le c_*N^{-3}$.
The ordering includes $s_2$ at the top, after decreasing $c_*$.
\end{corollary}
\begin{proof}
Subtract the two signs at each index in
\eqref{st-eq:polynomial-individual} and divide by $2\sqrt j$.
The remainder becomes $O(j|e|^{10})$. The common terms in the
normalized gaps cancel when taking the difference, as in the
proof of Theorem~\ref{st-thm:gap-variation}. This gives the
complex remainder estimate after composing with $e^2=S(I)$,
since $|I|\le c/k$ puts both indices in $|e|\le a_j/2$ and
$|e|\le a_k/2$ by the same lower bound on the amplitude scales.
The derivative estimate $-f'(x)\ge 233/(9800x^2)$ for $x\ge4$
and the boundary inequality
$\mathfrak f(3)-f(5)\ge(233/9800)(1/3-1/5)$ imply a leading
difference bounded below by a constant times $(k-j)/(jk)$.
The stated real interval makes the error smaller than half the
leading term. Since $(k-j)/(j k^2)\ge N^{-3}$ for $j<k\le N$,
the common interval follows.
At the fixed index $2$, the comparison's sixth-order expansion
and \eqref{st-eq:boundary} give the additional $(2/21)e^6$ term
in the normalized gap. The eighth-order normalized coefficients
at $k\ge3$ are bounded uniformly, and the remainder above is
$O(k|e|^{10})$. Consequently $s_2-s_k=(16/21)I^3+O(I^4+kI^5)$
on $0\le I\le c/k$ after decreasing $c$. This orders $s_2$ above
all those gaps on the stated finite-cutoff interval.
\end{proof}

\paragraph{Use of artificial intelligence.}
During the preparation of this work, the authors used ChatGPT (OpenAI) for editorial revision, reference checking, and supplementary computational verification. The authors reviewed and edited the output as needed and take full responsibility for the content of the published article.

\end{document}